\documentclass[11pt]{amsart}
\usepackage{amsmath,amssymb,amsfonts,url,mathptmx}
\numberwithin{equation}{section}

\newtheorem{thm}{Theorem}[section]
\newtheorem{lem}{Lemma}[section]
\newtheorem{rem}{Remark}[section]
\newtheorem{prop}{Proposition}[section]
\newtheorem{cor}{Corollary}[section]
\newcommand{\hdot}{^\text{\r{}}\hspace{-.33cm}H}

\begin{document}
\title[Toda system]{Convergence rate, location and $\partial_z^2$ condition for fully bubbling solutions to SU(n+1) Toda Systems} 

\author{Chang-shou Lin}
\address{Department of Mathematics\\
        Taida Institute of Mathematical Sciences\\
        National Taiwan University\\
         Taipei 106, Taiwan } \email{cslin@math.ntu.edu.tw}

\author{Juncheng Wei}
\address{Department of Mathematics\\
        Chinese University of Hong Kong\\
        Shatin, Hong Kong} \email{wei@math.cuhk.edu.hk}

\author{Lei Zhang}
\address{Department of Mathematics\\
        University of Florida\\
        358 Little Hall P.O.Box 118105\\
        Gainesville FL 32611-8105}
\email{leizhang@ufl.edu}

\date{\today}

\begin{abstract}
It is well known that the study of $SU(n+1)$ Toda systems is important not only to Chern-Simons models in Physics, but also to the understanding of holomorphic curves, harmonic sequences or harmonic maps from Riemann surfaces to $\mathbb C\mathbb P^n$. One major goal in the study of $SU(n+1)$ Toda system on Riemann surfaces is to completely understand the asymptotic behavior of fully bubbling solutions.
In this article we use a unified approach to study fully bubbling solutions to general $SU(n+1)$ Toda systems and we prove three major sharp estimates important for constructing bubbling solutions:  the closeness of blowup solutions to entire solutions, the location of blowup points and a $\partial_z^2$ condition.
\end{abstract}

\maketitle

\section{Introduction}
Let $(M,g)$ be a compact Riemann surface, in this article we consider the the following $SU(n+1)$ Toda system defined on $M$:
\begin{equation}\label{toda-1}
\Delta_g v_i+\sum_{j=1}^n a_{ij} H_j e^{v_j} -K(x)=4\pi \sum_{m} \gamma_{mi}\delta_{q_m},\quad 1\le i\le n
\end{equation}
where $\Delta_g$ is the Laplace-Beltrami operator ($-\Delta_g\ge 0$), $H_1,..H_n$ are positive smooth functions, $K$ is the Gauss curvature, $\delta_{q_m}$ stands for the Dirac measure at $q_m$, $A=(a_{ij})_{n\times n}$ is the following Cartan matrix:
$$A=\left(\begin{array}{ccccc}
2 & -1 & 0 & ... & 0 \\
-1& 2 & -1 & ... & 0 \\
0 & -1 & 2 & ... & 0 \\
\vdots & \vdots & ... & \vdots &  \\
0 & ... & -1 & 2 & -1 \\
0 & ... &   & -1 & 2
\end{array}
\right ) .
$$
The $SU(n+1)$ Toda system is well known to have close ties with many fields in Physics and Geometry. In Geometry the solutions of the Toda system are closely related to holomorphic curves (or harmonic sequences) of $M$ into $\mathbb C\mathbb P^n$. In the special case $M=\mathbb S^2$, the space of holomorphic curves of $\mathbb S^2$ to $\mathbb C\mathbb P^n$ is identical to the space of solutions to the $SU(n+1)$ system. The $q_m$s on the right hand side of (\ref{toda-1}) are ramification points of the corresponding holomorphic curve and $\gamma_{im}$ is the total ramificated index at $q_i$. See \cite{bolton,chern,griffith,guest} and the reference therein for discussions in detail.

On other hand in Physics, the analytic aspects of (\ref{toda-1}) are crucial for the understanding of the relativistic version of the non-abelian Chern-Simons models, see \cite{pacard-1,bartolucci-1,chen-lin-1,chen-lin-2,del-pino,li-cmp,lin-yan,nolasco-1,malchiodi-1,suzuki-1,tarantello-1} and the reference therein.

Using the Green's function
\begin{equation}\label{greenf}
\left\{\begin{array}{ll}
-\Delta_g G(y,\cdot)=\delta_p-1, \\
\int_M G(y,\eta)dV_g(\eta)=0.
\end{array}
\right.
\end{equation}
and a standard transformation (see \cite{lin-wei-zhao}) we can eliminate the singularity on the right hand side of (\ref{toda-1}) and rewrite (\ref{toda-1}) as
\begin{equation}\label{toda-2}
\Delta_g u_i+\sum_{j=1}^n a_{ij}\rho_j(\frac{h_j e^{u_j}}{\int_M h_j e^{u_j}dV_g}-1)=0, \quad \mbox{ in }\quad M,\quad i=1,..,n
\end{equation}
where $\rho_i>0$ are constants, $h_1,..,h_n$ are nonnegative continuous functions on $M$, and for convenience, we assume the volume of $M$ is $1$. It is easy to see that if $u=(u_1,..,u_n)$ is a solution, so is $(u_1+c_1,...,u_n+c_n)$ for arbitrary constants $c_1,..,c_n$, thus it is natural to use the following space for solutions to (\ref{toda-2}):
 $$\hdot (M)=\{u=(u_1,..,u_n)|\quad u_i\in H(M); \quad \int_M u_idV_g=0, \, i=1,..,n \}. $$

 System (\ref{toda-2}) is variational, as one can check immediately that it is the Euler-Lagrange equation of
 $$\Phi_{\rho}(u)=\frac 12\int_M\sum_{i,j=1}^n a^{ij}\nabla_g u_i\nabla_g u_j-\sum_{i=1}^n \rho_i \log \int_M h_i e^{u_i} dV_g $$
 where $(a^{ij})_{n\times n}$ is the inverse matrix of $A$.

 If (\ref{toda-2}) has only one equation, it is reduced to the following mean-field equation
 \begin{equation}\label{mean-f}
 \Delta_g u+\rho(\frac{he^u}{\int_M he^u dV_g}-1)=0, \quad \mbox{ in }\quad M.
 \end{equation}
 which has been extensively studied ( and fairly well understood) in the pass three decades because of its close connections with conformal geometry and Abelian Chern-Simons theory. See \cite{pacard-1,bartolucci-1,chen-lin-1,chen-lin-2,del-pino,lin-wei-03,li-cmp,shafrir,lin-yan,nolasco-1,malchiodi-1,suzuki-1,tarantello-1} and the reference therein for related discussions.

 In spite of its importance in theory and the profusion of applications, $SU(n+1)$ Toda system is well known for its analytical difficulties, since many fundamental tools, such as maximum principles, Pohozaev identities, that are very useful for single equations cannot be applied to Toda systems. Moreover, the solutions of $SU(n+1)$ Toda system have no symmetry whatsoever and are involved with too many parameters, for example, even for $SU(3)$ it takes 8 parameters to describe all the entire solutions in $\mathbb R^2$. When blowup solutions are considered, the asymptotic behavior of blowup solutions near an isolated blowup point is extremely complicated and the complexity increases significantly as the number of equations increases.

 The purpose of this article is to prove three major sharp estimates for fully bubbling solution (see (\ref{bas}) below) to general $SU(n+1)$ Toda systems: 1. All fully bubbling solutions are approximated by a sequence of global solutions with sharp error; 2. The gradient of certain functions must vanish sufficiently fast at blowup points; 3. There is a $\partial_z^2$ condition uniquely possessed by $SU(n+1)$ Toda systems. All these estimates are crucial for the understanding of bubble interactions and the construction of bubbling solutions in the future.

One of the major goals for (\ref{toda-2}) is first to determine the set of critical parameters when the blowup phenomenon occurs, and then derive a degree counting formula that depends on the topology of $M$.  The project has been successfully carried out for single Liouville equations ( \cite{chen-lin-1,chen-lin-2}) and Liouville systems (\cite{lin-zhang-1,lin-zhang-cpam,lin-zhang-jfa}).
Recently Lin-Wei-Yang \cite{lwya} derived a degree counting formula for $SU(3)$ Toda system. For general $SU(n+1)$-Toda systems there has been no degree counting formula and most of the progress so far is made on the $SU(3)$ Toda system (which has only two equations). For example, Jost-Lin-Wang \cite{jost-1} classified the limits of energy concentration for regular $SU(3)$ Toda systems, Lin-Wei-Zhang \cite{lin-wei-zhang} proved similar results for singular $SU(3)$ Toda systems, see \cite{malchiodi-1} and \cite{lwz-jems} for related discussions. In \cite{lin-wei-zhao} Lin, Wei and Zhao made significant progress by deriving the aforementioned three major sharp estimates for fully bubbling solutions to the $SU(3)$ Toda system, so in this article we removed the restriction on the number of equations and use a unified approach to extend Lin-Wei-Zhao's estimates to the general case.

Let $u^k=(u_1^k,...,u_n^k)$ satisfy
\begin{equation}\label{main-uk}
\Delta_g u_i^k+\sum_{j=1}^n a_{ij}\rho_j^k(\frac{h_j e^{u_j^k}}{\int_Mh_je^{u_j^k}dV_g}-1)=0,\quad \mbox{ on } M,\quad i=1,...,n.
\end{equation}
where $A=(a_{ij})_{n\times n}$ is the Cartan matrix, $M$ is a Riemann surface whose volume is assumed to be $1$ for convenience,
\begin{equation}\label{ha}
\mbox{ $h_1,...,h_n$ are
 positive smooth functions on $M$. }
 \end{equation}

 Let $p_1, ..., p_L$ be distinct blowup points. Our major assumption is
\begin{equation}\label{bas}
\mbox{Each blowup point is a fully bubbling blowup point}
\end{equation}
which means $(u_1^k,...,u_n^k)$ converges to a $SU(n+1)$ Toda system after appropriate scaling according to its maximum near the blowup point.
Specifically, let
\begin{equation}\label{ep}
M_k=\max_{i}\max_{x\in M}(u_i^k(x)-\int_Mh_ie^{u_i^k}dV_g),\quad \mbox{and }\quad \epsilon_k=e^{-\frac 12M_k}.
\end{equation}
Suppose $p_t^k$ satisfies $\lim_{k\to \infty}p_t^k=p_t$ (for $t=1,...,L$) and
$$\max_i\max_{x\in B(p_t,\delta)}u_i^k=\max_i u_i^k(p_t^k)$$
for some $\delta>0$ independent of $k$, then let
$$v_{t,i}^k(y):=u_i^k(p_t^k+\epsilon_k y)+2\log \epsilon_k,\quad i=1...,n $$
converges in $C^2_{loc}(\mathbb R^2)$ to
$$\Delta U_i+\sum_j a_{ij}h_i(p_t)e^{U_j}=0,\quad \mbox{ in } \mathbb R^2, \quad i=1,...,n $$
and $\int_{\mathbb R^2}e^{U_i}<\infty$. By the classification theorem of Lin-Wei-Ye \cite{lin-wei-ye}, $(U_1,...,U_n)$ is represented by $n^2+2n$ parameters. In other words, a blowup sequence is called fully bubbling if no component is lost after the scaling and taking the limit.

Our first main result is on the closeness between $\rho_i^k$ and $4\pi i(n+1-i)L$ when $k\to \infty$:
\begin{thm}\label{rho-1}
Let $(u_1^k,...,u_n^k)\in \hdot $ be a sequence of solutions to (\ref{main-uk}). Suppose (\ref{bas}) holds and we let $\epsilon_k$ and $h_i$ be described by (\ref{ep}) and (\ref{ha}), respectively. Then
for $i=1,...,n$
\begin{align*}
&\rho_i^k-4\pi i(n+1-i)L\\
=&\sum_{t=1}^L c_{k,i,t} (\Delta \log h_i(p_t^k)+8\pi L -2K(p_t^k))\epsilon_k^2|\log \epsilon_k |+O(\epsilon_k^2).
\end{align*}
where $K$ is the Gauss curvature on $M$, $L$ is the number of blowup points, $0<c_1<c_{k,i,t}\le c_2<\infty$ for $i=1,..,n$, $t=1,..,L$, and all $k$.
\end{thm}

Near a local point we write
$$G(y,\eta)=-\frac 1{2\pi }\log |y-\eta|+\gamma(y,\eta). $$

The second main result is on the locations of blowup points:
\begin{thm}\label{loc-bu-1} Under the same assumptions of Theorem \ref{rho-1}, the following vanishing condition holds
 for each blowup point $p_t^k$ ($t=1,...,L$):
$$8\pi \nabla_1 \gamma(p_t^k,p_t^k)+8\pi\sum_{l\neq t}\nabla_1 G(p_t^k,p_l^k)+\nabla (\log h_i)(p_t^k)=O(\epsilon_k) $$
where $i=1,...,n$ and $\nabla_1$ stands for the differentiation with respect to the first component.
\end{thm}

\begin{thm}\label{loc-bu-2} Under the same assumptions of Theorem \ref{rho-1}, the following $2n-2$ identities hold:
For $l=2,...,n$,
\begin{align}\label{2ndf1}
&\bigg (\Delta (\log h_{l-1})(p_t^k)+8\pi L-2K(p_t^k)\bigg )T_{l-1,k}^t \\
&+\bigg (\Delta (\log h_{l})(p_t^k)+8\pi L-2K(p_t^k)\bigg )T_{l,k}^t \nonumber \\
&+(l-1)l\pi(\partial_{11}-\partial_{22})(\log h_{l}-\log h_{l-1})(p_t^k)=O(\epsilon_k) \nonumber
\end{align}
and
\begin{align}\label{2ndf2}
&\bigg (\Delta (\log h_{l-1})(p_t^k)+8\pi L-2K(p_t^k)\bigg )\tilde T_{l-1,k}^t \\
&+\bigg (\Delta (\log h_{l})(p_t^k)+8\pi L-2K(p_t^k)\bigg )\tilde T_{l,k}^t \nonumber \\
&+2(l-1)l\pi \,\, \partial_{12}\,(\log h_{l}-\log h_{l-1})(p_t^k)=O(\epsilon_k) \nonumber
\end{align}
where $T_{l,k}^t$ and $\tilde T_{l,k}^t$ are obtained by differentiating parameters in approximating global solutions. See (\ref{Tlk}) for detail.
\end{thm}

Here we briefly describe $T_{l,k}^t$ and $\tilde T_{l,k}^t$. Around each $p_t^k$ the fully bubbling solution $u^k$ can be accurately approximated by a sequence of global solutions $U^k_t=(U_{1,t}^k,...,U_{n,t}^k)$ after scaling. The definition of $U^k_t$ is involved with $n^2+2n$ families of parameters ($c_{ij,t}^k$ and $\lambda_t^k$) which all have finite limits. $T_{l,k}^t$ and $\tilde T_{l,k}^t$ are obtained by differentiating $C_{n+2-i,n-i,t}^k$. Even though $C_{n+2-i,n-i,t}^k$ depends on $k$, it is uniformly bounded with respect to $k$. See \cite{lwz-jems} for details.

Theorem \ref{rho-1} and Theorem \ref{loc-bu-1} were established by Chen-Lin \cite{chen-lin-1} if the $SU(n+1)$ Toda system is reduced to the mean field equation (\ref{mean-f}). The $\partial_z^2$ condition in Theorem \ref{loc-bu-2} was first discovered by Lin-Wei-Zhao in \cite{lin-wei-zhao} for the $SU(3)$Toda system. It is interesting to observe that this condition does not exist when $n=1$ because $2n-2=0$ in this case. The reader can see that the major theorems:
Theorem \ref{rho-1},Theorem \ref{loc-bu-1} and Theorem \ref{loc-bu-2}, are extensions of the work of Lin-Wei-Zhao \cite{lin-wei-zhao} from $SU(3)$ to
 $SU(n+1)$ Toda system. Even for the case $n=2$, the estimates in
 Theorem \ref{loc-bu-2} is stronger than the corresponding estimate in \cite{lin-wei-zhao}.

 The proof of main theorems in this article is somewhat similar to the argument in \cite{lin-wei-zhao}. However our approach is systematic and has a number of
 new ideas. First we shall use the result in \cite{lwz-jems} as the initial step in our approximation. Then we study the algebraic structure of global solutions of $SU(n+1)$ Toda systems. It is well known that global solutions of $SU(n+1)$ Toda system are described by $n^2+2n$ parameters (see \cite{lin-wei-ye}). By carefully analyzing the leading terms of global solutions, we obtain certain families of solutions to the linearized system that are useful in our estimates. Third, our approach for pointwise estimate of blowup solutions is significantly simpler than the proof in \cite{lin-wei-zhao}. In particular we shall prove a pointwise sharp estimate for locally defined solutions in section three. The most essential reason of our argument is the classification of global solutions of $SU(n+1)$ Toda system and the non-degeneracy of the linearized system ( proved in \cite{lin-wei-ye}), which make us overcome all the difficulties from the lack of maximum principle and lack of symmetry.

 The organization of this paper is as follows. In chapter two we analyze the leading terms for global solutions of $SU(n+1)$ Toda systems. This part is based on the classification result of Lin-Wei-Ye \cite{lin-wei-ye}.  In Chapter three we study locally defined Toda systems with finite boundary oscillation on the boundary. We prove that all fully bubbling solutions can be accurately approximated by a set of global solutions. In this section we also establish the vanishing rate of certain function at the blowup point and the local version of the $\partial^2_z$ condition. This section is self-contained and may be interesting in its own right (see similar estimates for single equations in \cite{zhang-ccm,zhang-cmp} and for Liouville systems in \cite{lin-zhang-1},
 \cite{lin-zhang-cpam} and \cite{lin-zhang-jfa}). Then in Chapter four we use the results in Chapter three to prove all the main theorems in the introduction.

\bigskip

\noindent\textit{Acknowledgments.}
Part of the paper was finished when the third author was visiting
Taida Institute of Mathematical Sciences (TIMS) in July 2013 and University of British Columbia in November 2013. He would like to thank both institutes for their warm hospitality and generous financial support.

\section{Leading terms for global solutions}

In this section we identity the leading terms of the global solutions of
\begin{equation}\label{gU1}
\left\{\begin{array}{ll}
\Delta U_i+\sum_{j=1}^n a_{ij}e^{U_j}=0,\quad \mbox{ in }\quad \mathbb R^2,\quad i=1,...,n \\
\\
\int_{\mathbb R^2} e^{U_i}<\infty,\quad i=1,..,n,
\end{array}
\right.
\end{equation}
where $A=(a_{ij})_{n\times n}$ is the Cartan Matrix.

Lin-Wei-Ye\cite{lin-wei-ye} proved the following important classification result:

\emph{Theorem A (Lin-Wei-Ye): Let $(U_1,...,U_n)$ solve (\ref{gU1}) and
$$U^i=\sum_{j=1}^n a^{ij}U_j,\quad i=1,..,n $$
where $A^{-1}=(a^{ij})$ is the inverse matrix of $A$. $(U^1,....,U^n)$ are represented by $n^2+2n$ parameters:
$$U^1=-\log (\lambda_0+\sum_{i=1}^n \lambda_i |P_i(z)|^2) $$
where $\lambda_m>0$ for $m=0,...,n$,
$$ \lambda_0..\lambda_n=2^{-n(n+1)}\Pi_{1\le i\le j\le n}(j-i+1)^{-2}, $$
and
$$P_i(z)=z^i+\sum_{j=0}^{i-1}c_{ij}z^j, \quad \mbox{ for } \quad i=1,..,n. $$
Other components can be derived as follows:
Let $f=e^{-U^1}$, so
$$f=\lambda_0+\sum_{i=1}^n \lambda_i |P_i(z)|^2, $$
\begin{equation}\label{euk8}
e^{-U^k}=2^{k(k-1)}det_k(f),\quad \mbox{ for }\quad 2\le k\le n,
\end{equation}
where $det_k(f)=det(f^{p,q})_{0\le p,q\le k-1}$ for $1\le k\le n+1$, $f^{p,q}=\partial_{\bar z}^q\partial_z^pf$.
Moreover
$$|U_i(y)+4\log (1+|y|)|\le C, \quad y\in \mathbb R^2, \quad i=1,...,n. $$ }

\medskip

Recall that
$$a^{ij}=\frac{j(n+1-i)}{n+1},\quad \forall n\ge i\ge j\ge 1. $$
Since all
$$U_i(x)=-4\log |x|+O(1),\quad |x|>1,\quad \mbox{ for } i=1,...,n, $$
we have
$$U^i(x)=\sum_{j=1}^n a^{ij}U_j(x)=-4(\sum_{j=1}^n a^{ij})\log |x|+O(1),\quad |x|>1. $$
Elementary computation gives
\begin{align*}
4\sum_{j=1}^na^{ij}&=4(\sum_{j=1}^i a^{ij}+\sum_{j=i+1}^n a^{ij})\\
=&4(\sum_{j=1}^i \frac{j(n+1-i)}{n+1}+\sum_{j=i+1}^n \frac{i(n+1-j)}{n+1})\\
=&4(\frac{n+1-i}{n+1} \cdot \frac{i(i+1)}2+\frac{i}{n+1}\cdot \frac{(n-i)(n-i+1)}2 )\\
=&2i(n+1-i).
\end{align*}
Thus $U^i(x)=-2i(n+1-i)\log |x|+O(1)$ for $i=1,..,n$ and $|x|>1$.

In general we can write $e^{-U^m}$ as the following form:
\begin{prop}\label{euks}
\begin{align}\label{euk}
e^{-U^m}=&2^{m(m-1)}\lambda_n...\lambda_{n+1-m} \cdot \\
& \left | \begin{array}{ccc}
P_{n,3}  & ... &  P_{n+1-m,3} \\
... & ... & ... \\
P_{n,3}^{(m-1)} & ... & P_{n+1-m,3}^{(m-1)}
\end{array}
\right |
\cdot  \left | \begin{array}{ccc}
\bar P_{n,3} & ... &  \bar P_{n,3}^{(m-1)} \\
... & ... & ... \\
\bar P_{n+1-m,3} & ... & \bar P_{n+1-m,3}^{(m-1)}
\end{array}
\right | \nonumber \\
& \quad +ar^{2m(n+1-m)-2}+O(r^{2m(n+1-m)-3}), \quad m=1,..,n\nonumber
\end{align}
where $r=|z|$,$a\in \mathbb R$, $|\cdot |$ means the determinant of a matrix, $P_{m,3}^{(j)}$ stands for the first three terms of $\partial_z^jP_m$,for example,
$$P_{n,3}^{(1)}=nz^{n-1}+(n-1)c_{n,n-1}z^{n-2}+ (n-2)c_{n,n-2}z^{n-3},$$
 $\bar P_{m,3}^{(j)}$ stands for the first three terms of $\partial_{\bar z}^j\bar P_m$.
\end{prop}

\noindent{\bf Proof of Proposition \ref{euks}:}
We first prove the following lemma very closely related to the proof of Proposition \ref{euks}.
\begin{lem}\label{det1} Let $F(m)$ be a function defined by $F(1)=1$,
$$
F(m)=
\left | \begin{array}{cccc}
1 & 1 & ... & 1 \\
n & n-1 & ... & n+1-m \\
n(n-1) & (n-1)(n-2) & ... & (n+1-m)(n-m) \\
... & ... & ... &  ... \\
... & ... & ... &  ...  \\
\Pi_{j=0}^{m-2}(n-j) & \Pi_{j=1}^{m-1}(n-j) & ... & \Pi_{j=m-1}^{2m-3}(n-j)
\end{array}
\right |
$$ for $m\ge 2$.
Then \begin{equation}\label{fk}
F(m)=(-1)^{m(m-1)/2}(m-1)!(m-2)!...2!1!0!,\quad m\ge 1.
\end{equation}
\end{lem}

\noindent{\bf Proof of Lemma \ref{det1}:}
It is easy to see that $F(2)=-1$, so we assume $m\ge 3$.
We observe that if $2m+3\ge n$ some entries are zero, but it does not affect the proof. First we subtract $n+3-2m$ times row $m-1$ from row $m$, then we subtract $n+4-2m$ times row $m-2$ from row $m-1$, keep this process until subtracting $n+1-m$ times row $1$ from row $2$, we arrive at
\begin{align*}
&F(m)\\
=&
\left | \begin{array}{ccccc}
1 & 1 & ... & 1 & 1 \\
m-1 & m-2 & ... & 1 & 0 \\
(m-1)n & (m-2)(n-1) & ... & n+2-m & 0 \\
... & ... & ... & ... & ... \\
... & ... & ... & ... & ...  \\
(m-1)\Pi_{j=0}^{m-3}(n-j) & (m-2)\Pi_{j=1}^{m-2}(n-j) & ... & \Pi_{j=m-2}^{2m-5}(n-j) & 0
\end{array}
\right |
\end{align*}
Then by expanding along the last column and  taking out common factors from columns, we arrive at
\begin{align*}
&F(m)\\
=&(-1)^{m+1}(m-1)!\left | \begin{array}{cccc}
1 & 1 & ... & 1 \\
n & n-1 & ... & n+2-m \\
... & ... & ... & ... \\
\Pi_{j=0}^{m-4}(n-j) & \Pi_{j=1}^{m-3}(n-j) & ... & \Pi_{j=m-2}^{2m-6}(n-j)\\
\Pi_{j=0}^{m-3}(n-j) & \Pi_{j=1}^{m-2}(n-j) & ... & \Pi_{j=m-2}^{2m-5}(n-j)
\end{array}
\right |\\
=&(-1)^{m+1}(m-1)!F(m-1).
\end{align*}
Then it is easy to see that (\ref{fk}) holds.

Lemma \ref{det1} is established. $\Box$

\begin{rem}\label{rem-2p1}
If instead of using induction in the evaluation of $F(m)$, we employ the similar row reduction method on smaller matrices, then the matrix is reduced into the following form:
\begin{align*}
&F(m)\\
=&\left |\begin{array}{cccccc}
1 & 1 & ... & 1 & 1 & 1 \\
m-1 & m-2 & ... & 2 & 1 & 0 \\
(m-1)(m-2) & (m-2)(m-3) & ... & 2! & 0 & 0 \\
... & ... & ... & ... & ... & ... \\
(m-1)(m-2).. 2 & (m-2)! & ... & 0 & 0 & 0 \\
(m-1)! & 0 & ... & 0 & 0 & 0
\end{array}
\right |
\end{align*}
If we use $q_{ij}$ to represent the entry in row $i$ and column $j$, then $q_{ij}=0$ if $j>m+1-i$, $q_{i,m+1-i}=(i-1)!$ Then obviously the matrix can further be reduced to $q_{i,n+1-i}=(i-1)!$ and $q_{ij}=0$ if $j\neq n+1-i$ by row operations.
\end{rem}

Recall that $e^{-U^m}$ is determined by (\ref{euk8}).
We write $e^{-U^m}$ as
$$e^{-U^m}=2^{m(m-1)} \cdot det \bigg (( B_{\lambda}, C_{\lambda})\left(\begin{array}{c}
\bar B'\\
\bar C'
\end{array}
\right ) \bigg )$$
where
$$B_{\lambda}=\left ( \begin{array}{ccc}
\lambda_n P_n  & ... & \lambda_{n+1-m} P_{n+1-m} \\
... & ... & ... \\
\lambda_n P_n^{(m-1)} & ... & \lambda_{n+1-m} P_{n+1-m}^{(m-1)}
\end{array}
\right ) ,
$$
$$ C_{\lambda}=\left ( \begin{array}{ccc}
\lambda_{n-m} P_{n-m}  & ... &  \lambda_1 P_{1} \\
... & ... & ... \\
\lambda_{n-m} P_{n-m}^{(m-1)} & ... & \lambda_1 P_1^{(m-1)}
\end{array}
\right ) ,
$$
$$\bar B'=\left (\begin{array}{ccc}
\bar P_n & ... & \bar P_n^{(m-1)} \\
... & ... & ... \\
\bar P_{n+1-m} & ... & \bar P_{n+1-m}^{(m-1)}
\end{array}
\right ) ,
$$
and
$$\bar C'=\left (\begin{array}{ccc}
\bar P_{n-m} & ... & \bar P_{n-m}^{(m-1)} \\
... & ... & ... \\
\bar P_{n-m}^{(m-1)} & ... & \bar P_1^{(m-1)}
\end{array}
\right )
$$
For the matrix $( B_{\lambda}, C_{\lambda})$ we employ the same operations that we have employed in the proof of Lemma \ref{det1}. For each entry in this matrix, we just consider the leading term. For example, the first entry in row one is represented as $\lambda_n z^n(1+O(1/r))$ ($r=|z|$)(for convenience we shall use $1^*=1+O(1/r)$ in the computation below.  After doing the first set of row operations on $B_{\lambda}$ as in the proof of Lemma \ref{det1} and Remark \ref{rem-2p1} we can change $B_{\lambda}$ to
$$
\left (\begin{array}{ccccc}
\lambda_n z^n 1^* & \lambda_{n-1} z^{n-1}1^* & ... & ... & 0!\lambda_{n+1-m} z^{n+1-m}1^* \\
* & ... & ... & 1!\lambda_{n+2-m}z^{n+1-m}1^* & 0 \\
... & ... & ... & ... & ... \\
* & (m-2)! \lambda_{n-1}z^{n+1-m}1^* & 0 & ... &  0 \\
(m-1)! \lambda_n z^{n+1-m} 1^*& 0 & 0 &... & 0
\end{array}
\right )
$$
Note that $1^*$ changes from entry to entry since it just represents a quantity of the magnitude $1+O(1/r)$. Then we employ the second set of row operations to make all entries equal $0$ unless the the sum of the row number and the column number of the entry is $m+1$. Thus $B_{\lambda}$ is reduced to the following form: For some invertible matrix $Q$,
\begin{align*}
&QB_{\lambda}\\
=&
\left (\begin{array}{cccc}
0 & 0  & ... & 0!\lambda_{n+1-m} z^{n+1-m}1^* \\
0 & ...  & 1!\lambda_{n+2-m}z^{n+1-m}1^* & 0 \\
... & ...  & ... & ... \\
0 & (m-2)! \lambda_{n-1}z^{n+1-m}1^* &  ... &  0 \\
(m-1)! \lambda_n z^{n+1-m} 1^*& 0 & ... & 0
\end{array}
\right )  \\
\end{align*}
If the same columns operations are employed on $\bar B'$ we see that $\bar B'$ can be changed into this form:
$$
\bar B' \bar Q'
=\left (\begin{array}{ccccc}
0 & * & ... & * & (m-1)! \bar z^{n+1-m}1^* \\
0 & * & ... & (m-2)! \bar z^{n+1-m}1^* & 0 \\
... & ... & ... & ... & ...\\
0 & 1!\bar z^{n+1-m}1^* & ... & ... & 0 \\
0!\bar z^{n+1-m}1^* & 0 & ... & ... & 0
\end{array}
\right )
$$
When the same sets of row operations are employed on $C_{\lambda}$ we see that the entry in row 1 and column 1 has the highest $z$ power: $O(r^{n-m})$. All other entries have lower powers.
Clearly
\begin{align*}
&QB_{\lambda}\bar B'\bar Q'\\
=&\left(\begin{array}{cccc}
(0!)^2\lambda_{n+1-m}r^{2n+2-2m}1^* & 0 & ... & 0 \\
0 & (1!)^2\lambda_{n+2-m}r^{2n+2-2m}1^*  & ... & 0 \\
... & ... & ... & ... \\
0 & 0 & ... & ((m-1)!)^2\lambda_n r^{2n+2-2m}1^*
\end{array}
\right )
\end{align*}

Thus if we take $z^{n+1-m}$ from each row of $Q(B_{\lambda}, C_{\lambda})$ and $\bar z^{n+1-m}$ from each column of $\left(\begin{array}{c}
\bar B'\\
\bar C'
\end{array}
\right )\bar Q'
$, we
see that the leading term in $C_{\lambda}\bar C'$ is the $(1,1)$ entry of the order $O(r^{2n-2m})$ and it is of the form: $ar^{2n-2m}(1+O(1/r))$ for some $a\in \mathbb R$.  In the computation of $det(B_{\lambda}\bar B'+C_{\lambda}\bar C')$ only the $ar^{2n-2m}$ contributes to the term of order $O(r^{2m(n+1-m)-2})$. Then it is easy to see that each entry in $B_{\lambda}\bar B'$ only has at most three terms in its expansion that contribute to the form stated in (\ref{euk}).
Thus Proposition \ref{euks} is established. $\Box$

\medskip

The following proposition determines the second term in the expansion.
\begin{prop}\label{prop1} For $m\ge 2$
\begin{align*}
& \left | \begin{array}{cccc}
P_{n,3} & ... & ... & P_{n+1-m,3} \\
... & ... & ... & ... \\
P_{n,3}^{(m-1)} & ... & ... & P^{(m-1)}_{n+1-m,3}
\end{array}
\right |
\\
=&(-1)^{m(m-1)/2}(m-1)!...1! z^{m(n+1-m)} \\
&+\bigg (1+mc_{n+1-m,n-m}z^{-1}\bigg )
+O(|z|^{m(n+1-m)-2})
\end{align*}
\end{prop}
\begin{rem}
Proposition \ref{prop1} implies that for $m\ge 2$
\begin{align*}
& \left | \begin{array}{cccc}
\bar P_{n,3} & ... & ... & \bar P_{n+1-m,3} \\
... & ... & ... & ... \\
\bar P_{n,3}^{(m-1)} & ... & ... & \bar P^{(m-1)}_{n+1-m,3}
\end{array}
\right |
\\
=&(-1)^{m(m-1)/2}(m-1)!...1! \bar z^{m(n+1-m)}\\
&\bigg (1+m\bar c_{n+1-m,n-m}\bar z^{-1}\bigg )
+O(|z|^{m(n+1-m)-2}).
\end{align*}
\end{rem}

\noindent{\bf Proof of Proposition \ref{prop1}:}

It is easy to see that the leading term is of the order $O(r^{m(n+1-m)})$. Thus the second term is of the order $O(r^{m(n+1-m)-1})$. Then we see that the second term is a linear combination of $c_{n,n-1},...,c_{n+1-m,n-m})$:
$$\mbox{ The second term }=(\sum_{j=1}^m a_j c_{n+1-j,n-j} )z^{m(n+1-m)-1} $$
where $a_1,...,a_m$ are real numbers. We shall prove that only $a_m\neq 0$.
Now we compute the determinant of the first matrix in (\ref{euk}):
\begin{equation}\label{m1}
\mathbb M_P=\left | \begin{array}{ccc}
P_{n,3}  & ... &  P_{n+1-m,3} \\
... & ... & ... \\
P_{n,3}^{(m-1)} & ... & P_{n+1-m,3}^{(m-1)}
\end{array}
\right | .
\end{equation}
The leading term of $\mathbb M_P$ comes from the following determinant:
\begin{equation}\label{leadt}
\left | \begin{array}{cccc}
z^n & ... & z^{n+2-m} & z^{n+1-m} \\
nz^{n-1} & ... & (n+2-m)z^{n+1-m} & (n+1-m) z^{n-m} \\
... & ... & ... & ... \\
\Pi_{j=0}^{m-2}(n-j)z^{n+1-m} & ... & \Pi_{j=m-2}^{2m-4}(n-j)z^{n+3-2m} & \Pi_{j=m-1}^{2m-3}(n-j)z^{n+2-2m}
\end{array}
\right |
\end{equation}

After taking the smallest power of $z$ from each row and then taking out the common power of $z$ from each column we see that the determinant is equal to
\begin{align*}
&z^{m(n+1-m)}\\
&\left | \begin{array}{cccc}
1 & ... & 1 & 1 \\
n & ... & (n+2-m) & (n-m+1) \\
... & ... & ... & ... \\
\Pi_{j=0}^{m-2}(n-j)  & ... & \Pi_{j=m-2}^{2m-4}(n-j) & \Pi_{j=m-1}^{2m-3}(n-j)
\end{array}
\right |  \\
&=z^{m(n+1-m)}F(m)=z^{m(n+1-m)}(-1)^{m(m-1)/2}(m-1)!...1!.
\end{align*}

Next we compute the coefficient of $c_{n+1-j,n-j}$ in the $O(r^{m(n+1-m)-1})$ term. We see immediately that two columns of the corresponding matrix are the same if $j<m$. For example, in the computation of $c_{n,n-1}$, the first two columns of the corresponding matrix are the same. Thus only the coefficient of $c_{n+1-m,n-m}$ may be nonzero.

Since $c_{n+1-m,n-m}$ only comes from the last column of $\mathbb M_P$, we only need to consider the determinant of
$$
\left( \begin{array}{cccc}
z^n & ... & z^{n+2-m} & z^{n+1-m}+cz^{n-m} \\
nz^{n-1} & ... & (n+2-m)z^{n+1-m} & (n+1-m) z^{n-m}+c(n-m)z^{n-m-1} \\
... & ... & ... & ... \\
\Pi_{j=0}^{m-2}(n-j)z^{n+1-m} & ... & \Pi_{j=m-2}^{2m-4}(n-j)z^{n+3-2m} & \Pi_{j=m-1}^{2m-3}((n-j)z^{n+2-2m}
+c(n-j-1)z^{n+1-2m})
\end{array}
\right )
$$
where $c=c_{n+1-m,n-m}$.
Thus we need to compute the following determinant:
\begin{align} \label{cnk}
&c_{n+1-m,n-m}\\
&\left| \begin{array}{cccc}
z^n & ... & z^{n+2-m} & z^{n-m} \\
nz^{n-1} & ... & (n+2-m)z^{n+1-m} & (n-m) z^{n-m-1} \\
... & ... & ... & ... \\
\Pi_{j=0}^{m-2}(n-j)z^{n+1-m} & ... & \Pi_{j=m-2}^{2m-4}(n-j)z^{n+3-2m} & \Pi_{j=m}^{2m-2}(n-j)z^{n+1-2m}
\end{array}
\right |  \nonumber
\end{align}

Note that the powers of $z$ decrease by $1$ from column 1 to column $m-1$. However from column $m-1$ to column $m$, the powers decrease by $2$ in each row.
After taking out the power of $z$ just like what we did for (\ref{leadt}) the determinant is equal to
\begin{align*}
&c_{n+1-m,n-m} z^{m(n+1-m)-1}\\
&\left | \begin{array}{cccc}
1 & ... & 1 & 1 \\
n & ... & (n+2-m) & (n-m) \\
... & ... & ... & ... \\
\Pi_{j=0}^{m-2}(n-j) & ... & \Pi_{j=m-2}^{2m-4}(n-j) & \Pi_{j=m}^{2m-2}(n-j)
\end{array}
\right | \\
=&c_{n+1-m,n-m}G(m)z^{m(n+1-m)-1}.
\end{align*}
where
$$G(m)=\left | \begin{array}{cccc}
1 & ... & 1 & 1 \\
n & ... & (n+2-m) & (n-m) \\
... & ... & ... & ... \\
\Pi_{j=0}^{m-2}(n-j) & ... & \Pi_{j=m-2}^{2m-4}(n-j) & \Pi_{j=m}^{2m-2}(n-j)
\end{array}
\right |.
$$
To evaluate $G(m)$ we apply similar row reductions as we did on $F(m)$:
first subtract $n+2-2m$ times row $m-1$ from row $m$, then subtract $n+3-2m$ times row $m-2$ from row $m-1$, keep this process until subtracting $n-m$ times row $1$ from row two. Then we have
\begin{equation}\label{gk}
G(m)=(-1)^{m+1} m! F(m-1)=(-1)^{m(m-1)/2}m!(m-2)!...0!,\quad m\ge 2.
\end{equation}

Proposition \ref{prop1} is established. $\Box$

\medskip

Let
$$c_{n+1-m,n-m}=\alpha_m+\sqrt{-1}\beta_m, $$
by Proposition \ref{prop1} we write $e^{-U^m}$ as ( for $r=|z|$ large)
\begin{align}\label{eum}
&e^{-U^{m}}\\
=&2^{m(m-1)}\lambda_n...\lambda_{n+1-m}((m-1)!...0!)^2r^{2m(n+1-m)}\nonumber\\
&\bigg (1+2m(\alpha_m\cos\theta+\beta_m\sin\theta)r^{-1}+O(r^{-2})\bigg )\nonumber
\end{align}

\begin{rem} It is possible to write the form of $U_i$ by elementary computation. We list the result of some elementary computations:
$$b_{i,1}=\sum_{j=1}^n a_{ij} j(j-1)=\left\{\begin{array}{ll}
-2, \quad 1\le i\le n-1, \\
(n-1)(n+2),\quad i=n.
\end{array}
\right.
$$
$$b_{i,2}=\sum_j a_{ij} (\log \lambda_n+...+\log \lambda_{n+1-j})=\left\{\begin{array}{ll}\log \frac{\lambda_{n+1-i}}{\lambda_{n-i}},\quad
i<n, \\
\log \frac{\lambda_n...\lambda_2}{\lambda_1},\quad i=n.
\end{array}
\right.
$$
$$b_{i,3}=\sum_j a_{ij}(\log (0!)+...+\log (l-1)!)=\left\{\begin{array}{ll}
-\log i,\quad i<n, \\
\log (1!...(n-2)!((n-1)!)^2),\quad i=n.
\end{array}
\right.
$$
Then the expression of $U_i$ can be written as
$$-U_i(y)=b_{i,1}\log 2+b_{i,2}+2b_{i,3}+4\log |y|+O(1/r),\quad r=|y|>1, \quad i=1,..,n. $$
\end{rem}

Thus by (\ref{eum}) we have
\begin{align}\label{1-der}
-\frac{\partial U^{i}}{\partial \alpha_j}=\frac{2i}{r}\cos\theta \delta_{ij}+O(r^{-2})\\
-\frac{\partial U^{i}}{\partial \beta_j}=\frac{2i}{r}\sin\theta \delta_{ij}+O(r^{-2}),\quad \nonumber
\end{align}

By Proposition \ref{prop1}, in order to compute the $O(r^{2m(n+1-m)-2})$ term of $e^{-U^m}$ we first need to compute the $O(r^{m(n+1-m)-2})$ term of
$\mathbb M_p$ (see (\ref{m1})).
It is easy to observe that the $O(r^{m(n+1-m)-2})$ term is a linear combination of $c_{n-i,n-i-2}$ (for $i=0,...,n-2$) and $c_{n+1-i,n-i}c_{n+1-j,n-j}$ (where $i,j=1,...,n$). We focus on the coefficient of $c_{n-i,n-i-2}$. For $e^{-U^1}$ we only consider $c_{n,n-2}$ and $\bar c_{n,n-2}$. For $e^{-U^n}$ we only consider $c_{2,0}$ and $\bar c_{2,0}$. For $m=2,..,n-1$ we claim that the coefficient of $c_{n-i,n-i-2}$ and $\bar c_{n-i,n-i-2}$ is $0$ if $i<m-2$. Indeed, the coefficient of
$c_{n-i,n-i-2}$ comes from the determinant of the following matrix (we use $c$ to represent $c_{n-i,n-i-2}$)
\begin{align*}
&M_{c_{n+2-m,n-m}}\\
=&\left (\begin{array}{ccccc}
z^n & ... &  cz^{n-i-2} & ... & z^{n+1-m} \\
nz^{n-1} & ... & c(n-i-2)z^{n-i-3} & ... & (n+1-m)z^{n-m} \\
... & ... & ... & ... & ... \\
\Pi_{j=0}^{m-2}(n-j)z^{n+1-m} & ... & c\Pi_{j=i}^{i+m-2}(n-j)z^{n-i-1-m} & ... & \sum_{j=m-1}^{2m-3}(n-j)z^{n+2-2m}
\end{array}
\right )
\end{align*}
In $M_{c_{n+2-m,n-m}}$, the column with $c$ is the $i-th$ column, which comes from the second term of $z^{n-i}+cz^{n-i-2}$ and its derivatives. All other columns come from $z^j$ for $j\neq i$ and its derivatives. Then we see that if $i<m-2$, the $i-th$ column is a multiple of the $i+2$th column. Thus the determinant is $0$. Therefore we only need to determine the coefficient of $c_{n+2-m,n-m}$ and that of $c_{n+1-m,n-1-m}$ (and their conjugates) for $2\le m\le n-1$.

Before we consider the case $2\le m\le n$, we consider the case $m=1$ and $2$. By direct computation,
\begin{equation}\label{keq1}
-\frac{\partial U^{1}}{\partial \alpha_{2,2}}=\frac{2\cos (2\theta)}{r^2}+O(1/r^3), \quad
-\frac{\partial U^{1}}{\partial \beta_{2,2}}=\frac{2\sin (2\theta)}{r^2}+O(1/r^3).
\end{equation}
where $c_{n,n-2}=\alpha_{2,2}+\sqrt{-1}\beta_{2,2}$.

The case $m=2$ is involved with computations of $2\times 2$ matrices, so it is easy to verify from Proposition \ref{prop1} that
\begin{equation}\label{keq2}
-\frac{\partial U^2}{\partial \alpha_{2,2}}=\frac{-2\cos (2\theta)}{r^2}+O(r^{-3}),
\quad
-\frac{\partial U^2}{\partial \beta_{2,2}}=\frac{-2\sin (2\theta)}{r^2}+O(r^{-3}).
\end{equation}
and
\begin{equation}\label{keq2b}
-\frac{\partial U^2}{\partial \alpha_{3,2}}=\frac{6\cos (2\theta)}{r^2}+O(r^{-3}),
\quad
-\frac{\partial U^2}{\partial \beta_{3,2}}=\frac{6\sin (2\theta)}{r^2}+O(r^{-3}).
\end{equation}
where $c_{n-1,n-3}=\alpha_{3,2}+\sqrt{-1}\beta_{3,2}$.

\medskip

For $m\ge 3$, we take out the smallest power of $z$ from each row of $M_{c_{n+2-m,n-m}}$ and take the common power of $m$ from each column, then we
see that the term involving $c_{n+2-m,n-m}$ in $det(M_{c_{n+2-m,n-m}})$ is
\begin{align*}
&z^{m(n+1-m)-2}c_{n+2-m,n-m}\cdot \\
&\left |\begin{array}{ccccc}
1 & ... & 1 & 1 & 1 \\
n & ... & n+3-m & n-m & n+1-m \\
n(n-1) & ... & (n+3-m)(n+2-m) & (n-m)(n-m-1) & (n+1-m)(n-m)\\
... & ... & ... & ... & ... \\
\Pi_{j=0}^{m-2}(n-j) & ... & \Pi_{j=m-3}^{2m-5}(n-j) & \Pi_{j=m}^{2m-2}(n-j) & \Pi_{j=m-1}^{2m-3}(n-j)
\end{array}
\right |
\end{align*}
We apply the following row operations to the determinant above: Row $m$ minus $(n+2-2m)$ times row $m-1$, row $m-1$ minus $n+3-2m$ times row $m-2$,...,
row $2$ minus $n-m$ times row $1$. The penultimate column is reduced to $0$ except for the entry in row $1$. Then it is easy to see that starting from the second row, $m$ is a factor of all entries in the first column , $m-1$ is a common factor among all entries in the second column,..., $3$ is a common factor in column $m-2$. After taking out these common factors the determinant becomes
\begin{align*}
&z^{m(n+1-m)-2}c_{n+2-m,n-m}(-1)^m m(m-1)...3\cdot \\
&\left|\begin{array}{cccc}
1 & ... & 1 & 1 \\
n & ... & n+3-m & n+1-m \\
... & ... & ... & ... \\
\Pi_{j=0}^{m-3}(n-j) & ... & \Pi_{j=m-3}^{2m-6}(n-j) & \Pi_{j=m-1}^{2m-4}(n-j)
\end{array}
\right | \\
&=(-1)^m m...3 G(m-1)z^{m(n+1-m)-2}c_{n+2-m,n-m} \\
&=(-1)^{\frac{m(m-1)}2+1}m..3(m-1)!\Pi_{j=0}^{m-3}j! c_{n+2-m,n-m}z^{m(n+1-m)-2}.
\end{align*}

Correspondingly for $c_{n+1-m,n-1-m}$ ($m\le n-1$) we need to compute the following determinant:
$$\left |
\begin{array}{ccccc}
1 & ... & 1 & 1 & 1\\
n & ... & n+3-m & n+2-m & n-1-m \\
n(n-1) & ... & (n+3-m)(n+2-m) & (n+2-m)(n+1-m) & (n-1-m)(n-2-m) \\
... & ... & ... & ... & ... \\
\Pi_{j=0}^{m-2}(n-j) & ... & \Pi_{j=m-3}^{2m-5}(n-j) & \Pi_{j=m-2}^{2m-4}(n-j) & \Pi_{j=m+1}^{2m-1}(n-j)
\end{array}
\right |
$$
We employ similar row reductions as we have used for $M_{c_{n+2-m,n-m}}$: Row $m$ minus $n+1-2m$ times row $m-1$, ..., row $2$ minus $n-1-m$ times row $1$. Then the determinant is equal to
\begin{align*}
&\left |
\begin{array}{cccc}
1 & ... & 1 & 1 \\
(m+1) & ... & 3 & 0 \\
... & ... & ... & ... \\
(m+1)\Pi_{j=0}^{m-3}(n-j) & ... & 3\Pi_{j=m-2}^{2m-5}(n-j) & 0
\end{array}
\right | \\
=& (-1)^{m+1}(m+1)...3 F(m-1)\\
=&(-1)^{\frac{m(m-1)}2}(m+1)...3(m-2)!... 0!
\end{align*}
Let
$$
c_{n+2-m,n-m}=\alpha_{m,2}+\sqrt{-1}\beta_{m,2},\quad \mbox{ for }\quad  m=2,...,n, $$
then
$$ c_{n+1-m,n-1-m}=\alpha_{m+1,2}+\sqrt{-1}\beta_{m+1,2},\quad \mbox{ for } m=1,...,n-1. $$
Putting the estimates together we have, for $m\ge 3$,
the third term of $\mathbb M_p$ (see (\ref{m1})), which is of the order
$O(r^{m(n+1-m)-2})$, has two parts, one part is involved with products of $c_{n+i+1,n+i}$, which is not needed. The other part is
\begin{align*}
z^{m(n+1-m)-2}\bigg ((-1)^{\frac{m(m-1)}2+1}m...3 (m-1)!(m-3)!..0! c_{n+2-m,n-m} \\
+(-1)^{\frac{m(m-1)}2}(m+1)...3(m-2)!...0!c_{n+1-m,n-1-m}\bigg ).
\end{align*}
All other terms are error terms. Thus $e^{-U^m}$ can be written in the following form:
\begin{align}\label{euk2}
&e^{-U^m}\\
=&2^{m(m-1)}\lambda_n...\lambda_{n+1-m}((m-1)!...1!)^2 r^{2m(n+2-m)}\nonumber\\
&\bigg (1+2m(\alpha_m\cos \theta+\beta_m\sin\theta)/r   \nonumber \\
&+\big ((-m(m-1)\alpha_{m,2}+m(m+1)\alpha_{m+1,2})\cos (2\theta)\nonumber\\
&+(-m(m-1)\beta_{m,2}+m(m+1)\beta_{m+1,2})\sin (2\theta)\big )/r^2\nonumber \\
&+\mbox{ another } O(1/r^2) \mbox{ term } \bigg )\nonumber\\
&+O(r^{2m(n+2-m)-3}) \nonumber
\end{align}
The ``another $O(1/r^2)$ term" comes from the product of terms of the type $c_{n+1-i,n-i}$. (The reason that this term is not important is because when we differentiate global solutions with respect to $c_{n+1-i,n-i}$, this term is a high order error term when $r=|z|$ is large. )
Consequently for $3\le m\le n$,
\begin{align}\label{fre2}
-\frac{\partial U^m}{\partial \alpha_{m,2}}=-\frac{m(m-1)}{r^2}\cos (2\theta)+O(r^{-3}),\\
-\frac{\partial U^m}{\partial \beta_{m,2}}=-\frac{m(m-1)}{r^2}\sin (2\theta)+O(r^{-3}). \nonumber
\end{align}
Similarly for $3\le m\le n-1$ we have
\begin{align}\label{fre3}
-\frac{\partial U^m}{\partial \alpha_{m+1,2}}=\frac{m(m+1)}{r^2}\cos (2\theta)+O(r^{-3}),\\
-\frac{\partial U^m}{\partial \beta_{m+1,2}}=\frac{m(m+1)}{r^2}\sin (2\theta)+O(r^{-3}). \nonumber
\end{align}
Combining (\ref{keq1}),(\ref{keq2}),(\ref{keq2b}),(\ref{fre2}) and (\ref{fre3}) we have
\begin{align}\label{Udif}
-\frac{\partial U^m}{\partial \alpha_{m,2}}=-\frac{m(m-1)}{r^2}\cos (2\theta)+O(1/r^3),\quad m=2,...,n \\
-\frac{\partial U^m}{\partial \alpha_{m+1,2}}=\frac{m(m+1)}{r^2}\cos (2\theta)+O(1/r^3),\quad m=1,...,n-1 \nonumber \\
-\frac{\partial U^m}{\partial \beta_{m,2}}=-\frac{m(m-1)}{r^2}\sin (2\theta)+O(1/r^3),\quad m=2,...,n \nonumber \\
-\frac{\partial U^m}{\partial \beta_{m+1,2}}=\frac{m(m+1)}{r^2}\sin (2\theta)+O(1/r^3),\quad m=1,...,n-1 \nonumber
\end{align}

For $j=2,...,n$, let
$$\frac{\partial \bf U}{\partial \alpha_{j,2}}=(\frac{\partial U^1}{\partial \alpha_{j,2}},...,\frac{\partial U^n}{\partial \alpha_{j,2}})' $$
and
$$\frac{\partial \bf U}{\partial \beta_{j,2}}=(\frac{\partial U^1}{\partial \beta_{j,2}},...,\frac{\partial U^n}{\partial \beta_{j,2}})'. $$
By (\ref{Udif}) we have
\begin{equation}\label{ualpha}
-\frac{\partial \bf U}{\partial  \alpha_{j,2}}=\left(\begin{array}{c}
0 \\
.. \\
\frac{(j-1)j}{r^2}\cos(2\theta)\\
-\frac{(j-1)j}{r^2}\cos(2\theta)\\
0\\
...\\
0
\end{array}
\right)+O(r^{-3}),\quad j=2,...,n
\end{equation}
where the two non-zero entries are in rows $j-1$ and $j$, respectively. Similarly
\begin{equation}\label{ubeta}
-\frac{\partial \bf U}{\partial  \beta_{j,2}}=\left(\begin{array}{c}
0 \\
.. \\
\frac{(j-1)j}{r^2}\sin(2\theta)\\
-\frac{(j-1)j}{r^2}\sin(2\theta)\\
0\\
...\\
0
\end{array}
\right)+O(r^{-3}),\quad j=2,...,n
\end{equation}
where the two non-zero entries are in rows $j-1$ and $j$, respectively.

\section{Sharp estimates for locally defined Toda systems}

In this section we prove a pointwise, sharp estimate for locally defined Toda systems, which will be used to prove all the main results as easy consequences. Since this section has independent interest, we shall not choose notations different from other sections for simplicity.

Let
 $u^k=(u_1^k,...,u_n^k)$ be a sequence of solutions to
\begin{equation}\label{blowupu}
\left\{\begin{array}{ll}
\Delta u_i^k+\sum_{j=1}^n a_{ij} h_j e^{u_j^k}=0,\quad \mbox{ in }\Omega:=B(0,\sigma)\subset \mathbb R^2,\\
\\
\max_{x\in K}u_i^k(x)\le C(k), \quad \forall K\subset\subset \Omega\setminus \{0\},\quad i=1,...,n, \\
\\
|u_i^k(x)-u_i^k(y)|\le C_0,\quad \forall x,y\in \partial \Omega, \quad i=1,...,n \\
\\
\mbox{ There exists $C$ independent of $k$ such that } \int_{B_1}h_ie^{u_i^k}\le C, \quad i=1,..,n.
\end{array}
\right.
\end{equation}
where $A=(a_{ij})_{n\times n}$ is the Cartan matrix. $h_1,...,h_n$ are positive smooth functions in $\Omega$:
\begin{equation}\label{ah}
\frac 1C_1\le h_i(x)\le C_1,\quad \|\nabla h_i\|_{C^3(\Omega)}\le C_1.
\end{equation}
From the second equation of (\ref{blowupu}) we see that $0$ is the only possible blowup point. We write the equation for $u_i^k$ as
$$\Delta (u_i^k(x)+\log h_i(0))+\sum_j a_{ij}\frac{h_j(x)}{h_j(0)}e^{u_j(x)+\log h_j(0)}=0. $$
Our major assumption is that $u^k$ is a fully bubbling sequence: Let
$$v_i^k(y)=u_i^k(\epsilon_ky)+\log h_i(0)+2\log \epsilon_k,\quad i=1,...,n, \quad e^{-\frac 12 \epsilon_k}=\max_{i}\max_{\bar \Omega}(u_i^k+\log h_i(0)). $$
Suppose $\epsilon_k\to 0$ and $v^k=(v_1^k,...,v_n^k)$ converges in $C^2_{loc}(\mathbb R^2)$ to $U=(U_1,...,U_n)$, which satisfies (\ref{gU1}).

Let $\phi_i^k$ be defined by
$$\left\{\begin{array}{ll}
\Delta \phi_i^k=0,\quad \mbox{ in }\quad \Omega, \\
\phi_i^k(x)=u_i^k(x)-\frac{1}{2\pi \sigma}\int_{\partial \Omega}u_i^k dS \quad \mbox{ on }\partial \Omega, \quad i=1,...,n .
\end{array}
\right.
$$
Clearly $\phi_i^k(0)=0$, all derivatives of $\phi_i^k$ are uniformly bounded over all fixed compact subsets of $\Omega$ because
$u_i^k$ has bounded oscillation on $\partial B_1$.

Clearly $v^k=(v_1^k,...,v_n^k)$ satisfies
$$\Delta v_i^k(y)+\sum_j a_{ij}\frac{h_j(\epsilon_k y)}{h_j(0)}e^{v_j^k(y)}=0,\quad \mbox{ in }\quad \Omega_k:=\{y;\quad \epsilon_k y\in \Omega\}. $$
Let
$$\tilde v_i^k(\cdot )=v_i^k(\cdot)-\phi_i^k(\epsilon_k \cdot ), \quad \mbox{ in }\quad \Omega_k $$
 and
 $$h_i^k(\cdot)=\frac{h_i(\cdot)}{h_i(0)}e^{\phi_i^k(\cdot)},\quad x\in \Omega, $$
 then we have
 $$\Delta \tilde v_i^k(y)+\sum_{j=1}^n a_{ij} h_j^k(\epsilon_k y)e^{\tilde v_j^k(y)}=0, \quad \mbox{ in }\quad \Omega_k. $$
It is proved in \cite{lwz-jems} that there exists a sequence $U^k=(U_1^k,...,U_n^k)$ such that the following properties hold:
\begin{enumerate}
\item
$U^k=(U_1^k,...,U_n^k)$ satisfies (\ref{gU1}).
\item Let $c_{ij}^k$ and $\lambda_i^k$ be the $n^2+2n$ families of parameters in the definition of $U^k$, then they all converge to finite limits: $c_{ij}^k\to c_{ij}$ and
$\lambda_i^k\to \lambda_i$.
\item Let
\begin{equation}\label{U-up-ik}
U^{i,k}=\sum_{j=1}^n a^{ij}U_j^k \quad \mbox{and}\quad \tilde v^{i,k}=\sum_{j=1}^n a^{ij}\tilde v_j^k,
\end{equation}
there are $n^2+2n$ distinct points $p_1,...,p_{n^2+2n}\in \mathbb R^2$, all independent of $k$, such that
\begin{equation}\label{uagree}
U^{1,k}(p_l)=\tilde v^{1,k}(p_l),\quad l=1,...,n^2+2n.
\end{equation}
\item
There is a $C$ independent of $k$ such that
\begin{equation}\label{bad-ap}
|\tilde v_i^k(y)-U_i^k(y)|\le C(\delta)\epsilon_k (1+|y|),\quad y\in \Omega_k,\quad i=1,...,n.
\end{equation}
\end{enumerate}

The main result of the locally defined system is:

\begin{thm}\label{localthm}
Let $u^k$, $\phi^k$, $h_1,...h_n$, $\tilde v_i^k$, $U_i^k$ be described as above. There exists $C>0$ independent of $k$ such that
\begin{equation}\label{pw}
|\tilde v_i^k(y)-U_i^k(y)|\le C\epsilon_k^2(1+|y|),\quad \mbox{ for }\quad y\in \Omega_k.
\end{equation}
Moreover,
\begin{equation}\label{1van}
(\nabla \log h_i^k)(0)=\frac{\nabla h_i(0)}{h_i(0)}+\nabla \phi_i^k(0)=O(\epsilon_k), \quad i=1,...,n
\end{equation}
and for $l=2,...,n$, the following $2n-2$ identities hold:
\begin{align}\label{2van}
&\Delta (\log h_{l-1}^k)(0)\int_{\mathbb R^2}(-\frac{\partial U^{l-1,k}}{\partial \alpha^k_{l,2}})dx
+\Delta (\log h_l^k)(0)\int_{\mathbb R^2}(-\frac{\partial U^{l,k}}{\partial \alpha_{l,2}^k})dx\\
&+(l-1)l\pi(\partial_{11}-\partial_{22})(\log h_l^k-\log h_{l-1}^k)(0)=O(\epsilon_k) \nonumber \\
&\Delta (\log h^k_{l-1})(0)\int_{\mathbb R^2}(-\frac{\partial U^{l-1,k}}{\partial \beta^k_{l,2}})dx
+\Delta (\log h^k_{l})(0)\int_{\mathbb R^2}(-\frac{\partial U^{l,k}}{\partial \beta^k_{l,2}})dx \nonumber \\
&+2(l-1)l\pi\partial_{12}(\log h^k_{l}-\log h_{l-1}^k)(0)=O(\epsilon_k) \nonumber
\end{align}
where $U^{i,k}$ is defined in (\ref{U-up-ik}) and
$$c_{n+2-i,n-i}^k=\alpha^k_{i,2}+\sqrt{-1}\beta^k_{i,2} \quad \mbox{ for } i=2,...,n.
$$
\end{thm}

\noindent{\bf Proof of Theorem \ref{localthm}:}

\subsection{Vanishing theorem on the coefficient of the first frequency}

Let
$$w_i^k(y)=\tilde v_i^k(y)-U_i^k(y),\quad \mbox{ for } y\in \Omega_k, \quad i=1,...,n. $$
It is proved in \cite{lwz-jems} that
\begin{equation}\label{wik1}
|w_i^k(y)|\le C\epsilon_k (1+|y|).
\end{equation}
Using (\ref{wik1}) we write the equation for $w_i^k$ as
\begin{align}\label{wik2}
&\Delta w_i^k+\sum_j a_{ij} h_j^k(\epsilon_ky) e^{\xi_j^k}w_j^k(y)\\
=&-\epsilon_k \sum_j a_{ij} (\partial_1 h_j^k(0)y_1+\partial_2 h_j^k(0)y_2)e^{U_j^k}+O(\epsilon_k^2)(1+|y|)^{-2} \nonumber
\end{align}
where $e^{\xi_i^k}$ is obtained from mean value theorem.
Let
$$w^{i,k}=\sum_j a^{ij}w_j^k$$
where $(a^{ij})_{n\times n}=A^{-1}$. Then we have
\begin{equation}\label{wik-up}
\Delta w^{i,k}+e^{U_i^k}w_i^k=-\epsilon_k(\partial_1 h_i^k(0)y_1+\partial_2 h_i^k(0)y_2)e^{U_i^k}+O(\epsilon_k^2)(1+|y|)^{-2}.
\end{equation}
where we have used  $w_i^k=O(\epsilon_k)(1+|y|)$, $h_i^k(0)=1$, and
$$e^{\xi_i^k}-e^{U_i^k}=O(\epsilon_k)(1+|y|)^{-3}.$$
We shall multiply to both sides of (\ref{wik-up}) $\phi_i$, which solves
\begin{equation}\label{phi-1}
\Delta \phi_i+\sum_j a_{ij} e^{U_j^k}\phi_j=0,\quad \mbox{ in }\mathbb R^2.
\end{equation}
$\phi_i$ will be chosen to satisfy
\begin{equation}\label{phi-decay}
\phi_i(y)=O(1/|y|),\quad |y|>1, \quad |\nabla \phi_i(y)|=O(1/|y|^2),\quad |y|>1.
\end{equation}
Correspondingly we define
$$\phi^i=\sum_j a^{ij}\phi_j. $$
For simplicity we do not include $k$ in $\phi_i$ and $\phi^i$.
Then $\phi^i$ satisfies
\begin{equation}\label{phi-up}
\Delta \phi^i+e^{U_i^k}\phi_i=0,\quad i=1,...,n,\quad \mbox{ in }\quad \mathbb R^2.
\end{equation}

Here we recall that $\Omega_k=B(0,\epsilon_k^{-1}\sigma)$ and we define $\tilde \Omega_k$ as $B(0,\frac 12 \epsilon_k^{-1}\sigma)$. By multiplying $\phi_i$ to both sides of (\ref{wik-up})and integrating on $\tilde \Omega_k$, we obtain the following equation:
\begin{align}\label{tem-eq1}
&\sum_i \int_{\tilde \Omega_k} (\Delta w^{i,k}+e^{U_i^k}w_i^k)\phi_i\\
=&-\sum_i \epsilon_k \int_{\tilde \Omega_k}(\partial_1 h_i^k(0)y_1+\partial_2 h_i^k(0)y_2)e^{U_i^k}\phi_i+O(\epsilon_k^2) \nonumber
\end{align}
On the left hand side

\begin{align*}
&\sum_i\int_{\tilde \Omega_k}(\Delta w^{i,k}+e^{U_i^k}w_i^k)\phi_i \\
&=\sum_i\int_{\partial \tilde \Omega_k}(\partial_{\nu}w^{i,k}\phi_i-w^{i,k}\partial_{\nu}\phi_i)+\sum_i\int_{\tilde \Omega_k}
(\Delta \phi_iw^{i,k}+e^{U_i^k}w_i^k\phi_i).
\end{align*}
Since
$$\sum_i \Delta \phi_i w^{i,k}=\sum_i \Delta \phi_i(\sum_l a^{il}w_l^k), $$
and $a^{il}=a^{li}$, by interchanging $i$ and $l$ we have
$$\sum_i\int_{\tilde \Omega_k}\Delta \phi_i(\sum_l a^{il}w_l^k)=\sum_i\int_{\tilde \Omega_k}\Delta \phi^i w_i^k. $$
Thus
\begin{align}\label{IGB}
&\sum_i\int_{\tilde \Omega_k}(\Delta w^{i,k}+e^{U_i^k}w_i^k)\phi_i \\
=&\sum_i\int_{\partial \tilde \Omega_k}(\partial_{\nu}w^{i,k}\phi_i-w^{i,k}\partial_{\nu}\phi_i)
+\sum_i\int_{\tilde \Omega_k}(\Delta \phi^i+e^{U_i^k}\phi_i)w_i^k \nonumber \\
=&\sum_i\int_{\partial \tilde \Omega_k}(\partial_{\nu}w^{i,k}\phi_i-w^{i,k}\partial_{\nu}\phi_i). \nonumber
\end{align}

In order to evaluate different terms on $\partial \tilde \Omega_k$ we need the following estimate
\begin{align} \label{1dew}
\nabla w_i^k=O(\epsilon_k^2)\quad \mbox{ on }\partial \tilde \Omega_k \\
w_i^k=c_i^k+O(\epsilon_k), \,\, \mbox{ on } \partial \tilde \Omega_k,
\nonumber
\end{align}
where $c_i^k$ are uniformly bounded constants.

The proof of (\ref{1dew}) is based on the following Green's representation of $w_i^k$:
\begin{equation}\label{ntem-2}
w_i^k(y)=\int_{\Omega_k}G_k(y,\eta)(-\Delta w_i^k)-\int_{\partial \Omega_k}\partial_{\nu}G_k(y,\eta)w_i^k(\eta)dS_{\eta}
\end{equation}
where $G_k$ is the Green's function on $\Omega_k$:
$$G_k(y,\eta)=-\frac 1{2\pi}\log |y-\eta |+\frac 1{2\pi} \log (\frac{|y|}{\epsilon_k^{-1}}|\frac{\epsilon_k^{-2}y}{|y|^2}-\eta|),\quad y,\eta\in \Omega_k. $$
The second term in (\ref{ntem-2}) is a harmonic function on $\Omega_k$ with $O(\epsilon_k)$ oscillation on $\partial \Omega_k$. Thus its derivative is $O(\epsilon_k^2)$ on $\partial \tilde \Omega_k$. Thus we only need to prove the derivative of the first term of (\ref{ntem-2}) is $O(\epsilon_k^2)$ on $\partial \tilde \Omega_k$. To this end we just need to verify that
\begin{align*}
\int_{\tilde \Omega_k}\nabla_y G_k(y,\eta)\bigg (e^{U_i^k}w_i^k+\epsilon_k(\partial_1h_i^k(0)\eta_1+\partial_2h_i^k(0)\eta_2)e^{U_i^k}\\
+O(\epsilon_k^2)(1+|\eta |)^{-2}\bigg )d\eta=O(\epsilon_k^2),\quad y\in \partial \tilde \Omega_k.
\end{align*}
The desired estimate follows from (\ref{bad-ap}), $e^{U_i^k(\eta)}=O(|\eta |^{-4})$ and standard estimates. Thus (\ref{1dew}) is verified.

$\phi^i$ will be chosen to satisfy
\begin{equation}\label{phi-eqa}
\phi^i=(a_i^k\cos \theta+b_i^k \sin \theta)/r+O(r^{-2}).
\end{equation}
Then from (\ref{1dew}) and (\ref{phi-eqa}) we have
$$\sum_i\int_{\partial \tilde \Omega_k}(\partial_{\nu}w^{i,k}\phi_i-w^{i,k}\partial_{\nu}\phi_i)=O(\epsilon_k^2). $$

Now we focus on the right hand side of (\ref{tem-eq1}). For each $i$, by (\ref{phi-up})
$$e^{U_i^k}\phi_i=-\Delta \phi^i. $$
Thus integration by parts gives
\begin{align}\label{13nov6e1}
&\sum_i \int_{\tilde \Omega_k}(\partial_1h^k_i(0)y_1+\partial_2h^k_i(0)y_2)e^{U_i^k}\phi_i \\
=&\sum_i \int_{\partial \tilde \Omega_k}(y_1 \partial_1 h^k_i(0)+y_2 \partial_2 h^k_i(0))(\frac{\phi^i}{|y|}-\partial_{\nu}\phi^i). \nonumber
\end{align}
Using polar coordinates, what we need to show is
\begin{equation}\label{need1}
\sum_i\int_0^{2\pi}(\cos\theta \partial_1h^k_i(0)+\sin\theta \partial_2 h^k_i(0))(\phi^i-r\partial_r \phi^i)rd\theta=O(\epsilon_k)
\end{equation}
where $r=\frac{\sigma}2 \epsilon_k^{-1}$.

Recall that $(U^{1,k},...,U^{n,k})$ are described by parameters $\lambda_m^k$ ($m=0,...,n$) and $c_{ij}^k$ ($i>j$).
Let
$$c^k_{n+1-m,n-m}=\alpha^k_m+\sqrt{-1}\beta^k_m, $$
and
$$\Phi_{\alpha_j}=\left ( \begin{array}{c}
\phi^1_{\alpha_j} \\
... \\
... \\
\phi^n_{\alpha_j}
\end{array}
\right )
=
\left (\begin{array}{c}
-\frac{\partial U^{1,k}}{\partial \alpha^k_j} \\
... \\
... \\
-\frac{\partial U^{n,k}}{\partial \alpha^k_j}
\end{array}
\right )=
\left (\begin{array}{c}
2\delta_1^j\cos \theta/r\\
...\\
...\\
2n\delta_n^j\cos\theta/r
\end{array}
\right )+O(1/r^2),
$$

$$\Phi_{\beta_j}=\left ( \begin{array}{c}
\phi^1_{\beta_j} \\
... \\
... \\
\phi^n_{\beta_j}
\end{array}
\right )
=
\left (\begin{array}{c}
-\frac{\partial U^{1,k}}{\partial \beta^k_j} \\
... \\
... \\
-\frac{\partial U^{n,k}}{\partial \beta^k_j}
\end{array}
\right )=
\left (\begin{array}{c}
2\delta_1^j\sin \theta/r\\
...\\
...\\
2n\delta_n^j\sin\theta/r
\end{array}
\right )+O(1/r^2)
$$
where we have used (\ref{1-der}).
For $j=1,...,n$ we have $2n$ sets of solutions to the linearized system.

From (\ref{need1}) we observe that if
$$\phi^i=\frac 1r(d_i \cos \theta+ q_i \sin \theta)+O(1/r^2) $$
we have
\begin{align*}
&\sum_i \int_0^{2\pi} (\cos \theta \partial_i h^k_i(0)+\sin \theta \partial_2 h^k_i(0))r(\phi^i-r\partial_r \phi^i)d\theta\\
=& \pi \sum_i (d_i \partial_1 h^k_i(0)+q_i \partial_2 h^k_i(0))+O(\epsilon_k)
\end{align*}
where $r\sim \epsilon_k^{-1}$.

Replacing $\phi^i$ by $\Phi_{\alpha_j}$ and $\Phi_{\beta_j}$ for $j=1,...,n$, we have
\begin{equation}\label{first-freq}
\partial_1 h^k_i(0)=O(\epsilon_k), \quad \partial_2 h^k_i(0)=O(\epsilon_k),\quad i=1,...,n.
\end{equation}

\bigskip

\subsection{The vanishing rate for the second frequency terms}

The estimates of $w_i^k$ can be improved to the following form:
\begin{prop}\label{wik-ip} There exists $C>0$ independent of $k$ such that
$$|w_i^k(y)|\le C\epsilon_k^2(1+|y|),\quad y\in \Omega_k.$$
\end{prop}

\noindent{\bf Proof of Proposition \ref{wik-ip}:}

With the vanishing rate of $\nabla h^k_i(0)$.
The equation for $w_i^k$ can now be written as
$$\Delta w^{i,k}+e^{\xi_i^k}w_i^k=O(\epsilon_k^2)(1+|y|)^{-2},\quad y\in \Omega_k. $$
Recall that $\tilde v_i^k=constant$ on $\partial \Omega_k$,
$$U_i^k(y)=-4\log |y|+c_{i,k}+(a_i^k\cos \theta+b_i^k\sin\theta )/r+O(1/r^2),\quad \mbox{ on } \partial \Omega_k $$
for some uniformly bounded constants $c_{i,k}$, $a_i^k$, $b_i^k$. Let $\psi_i^k$ be a harmonic function in $\Omega_k$ that satisfies
$\tilde v_i^k-U_i^k-\psi_i^k=$  constant on $\partial \Omega_k$ and $\psi_i^k(0)=0$. Since the oscillation of $\tilde v_i^k-U_i^k$ is $O(\epsilon_k)$ on $\partial \Omega_k$ we have
\begin{equation}\label{psik}
|\psi_i^k(y)|\le C\epsilon_k^2|y|, \quad y\in \Omega_k.
\end{equation}

Set
$$\tilde w_i^{k}=v_i^k-U_i^k-\psi_i^k$$
Then $\tilde w_i^{k}=constant$ on $\partial \Omega_k$.
The equation for $\tilde w_i^{k}$ is
$$\Delta \tilde w^{i,k}+e^{\xi_i^k}\tilde w_i^{k}=O(\epsilon_k^2)(1+|y|)^{-2},\quad y\in \Omega_k. $$
where $\tilde w^{i,k}=\sum_j a^{ij}\tilde w_j^{k}$.
Given $\delta\in (0,1)$, we let
$$\Lambda_k=\sup_i\sup_{x\in \Omega_k}\frac{|\tilde w^{i,k}(x)|}{\epsilon_k^2(1+|y|)^{\delta}}. $$
Our goal is to prove that $\Lambda_k\le C$. Suppose this is not the case, we use $y_k$ to denote where the maximum is attained on $\bar
\Omega_k$. Let
\begin{equation}\label{hat-wik}
\hat w^{i,k}(y)=\frac{\tilde w^{i,k}(y)}{\Lambda_k \epsilon_k^2(1+|y_k|)^{\delta}}.
\end{equation}
It is easy to see from the definition of $\Lambda_k$ that
\begin{equation}\label{hat-wik-2}
|\hat w^{i,k}(y)|=|\frac{\hat w^{i,k}(y)}{\Lambda_k \epsilon_k^2(1+|y|)^{\delta}}\frac{(1+|y|)^{\delta}}{(1+|y_k|)^{\delta}}|\le \frac{(1+|y|)^{\delta}}{(1+|y_k|)^{\delta}}.
\end{equation}
The equation for $\hat w^{i,k}$ can be written as
$$-\Delta \hat w^{i,k}=\frac{O((1+|y|)^{-4+\delta})}{(1+|y_k|)^{\delta}}+\frac{O((1+|y|)^{-2})}{\Lambda_k (1+|y_k|)^{\delta}}. $$
First we claim
$
|y_k|\to \infty.
$
 By way of contradiction we assume that $y_k\to y^*\in \mathbb R^2$. In this case we see that along a subsequence
$\hat w_i^k$ converges to $(\phi_1,...,\phi_n)$ in $C^2_{loc}(\mathbb R^2)$ that satisfies
\begin{equation}\label{limit-phi}
\left\{\begin{array}{ll}
\Delta \phi_i+\sum_j a_{ij}e^{U_j}\phi_j=0,\quad \mathbb R^2 \\
|\phi_i(y)|\le C(1+|y|)^{\delta},\quad y\in \mathbb R^2, \\
\phi_1(p_l)=0,\quad \mbox{ for distinct } p_1,..., p_{n^2+2n}\in \mathbb R^2.
\end{array}
\right.
\end{equation}
Indeed, from (\ref{hat-wik-2}) the second equation of (\ref{limit-phi}) follows. The reason that the third equation of (\ref{limit-phi}) holds is because by (\ref{uagree}) we have $w^{1,k}(p_l)=0$ for $l=1,..,n^2+2n$. Then from the definition of $\hat w^{1,k}$ in (\ref{hat-wik}) and (\ref{psik}) we have
$$\hat w^{1,k}(p_l)=O(1/\Lambda_k),\quad l=1,...,n^2+2n. $$
Thus the third equation of (\ref{limit-phi}) holds.

By the classification theorem of Lin-Wei-Ye \cite{lin-wei-ye}, $\phi_i\equiv 0$ for all $i$. Thus a contradiction to $\phi_j(y^*)=\pm 1$ for some $j$. Therefore we only need to rule out the possibility that $y_k\to \infty$.

By Green's representation formula and $\hat w^{i,k}=constant$ on $\partial \Omega_k$, we have
$$\hat w^{i,k}(y_k)=\int_{\Omega_k}(G_k(y_k,\eta)-G_k(0,\eta))
(\frac{O((1+|\eta|)^{-4+\delta})}{(1+|y_k|)^{\delta}}+\frac{O((1+|\eta|)^{-2})}{\Lambda_k (1+|y_k|)^{\delta}})d\eta. $$
Using the standard estimates on the Green's function $G_k$ (see Lemma 3.2 of \cite{lin-zhang-jfa})
it is easy to see that the right hand side is $o(1)$, a contradiction to the fact that $\hat w^{i,k}(y_k)=\pm 1$ for some $i$. Proposition \ref{wik-ip} is established. $\Box$

\medskip

Using Proposition \ref{wik-ip} we write the equation for $w_i^k$ as
\begin{align*}
&\Delta w^{i,k}+e^{U_i^k}w_i^k \\
=&-\epsilon_k(\partial_1h^k_i(0)y_1+\partial_2h^k_i(0)y_2)e^{U_i^k}\\
&-(\partial_{11}h^k_i(0)y_1^2+2\partial_{12}h^k_i(0)y_1y_2+\partial_{22}h^k_i(0)y_2^2)e^{U_i^k}\epsilon_k^2+O(\epsilon_k^3)(1+|y|)^{-1}\\
=&-\epsilon_k(\partial_1h^k_i(0)y_1+\partial_2h^k_i(0)y_2)e^{U_i^k}\\
&-\bigg (\frac 12\Delta h^k_i(0)+\frac 12(\partial_{11}h^k_i(0)-\partial_{22}h^k_i(0))\cos(2\theta)+\partial_{12}h^k_i(0)\sin (2\theta)\bigg )r^2e^{U_i^k}\epsilon_k^2\\
&\quad +O(\epsilon_k^3)(1+|y|)^{-1}.
\end{align*}
Multiplying $\phi_i$ to both sides, taking the summation on $i$ and integrating on $\tilde \Omega_k$, we change the left hand side as before to
$$\sum_i \int_{\partial \tilde \Omega_k} (\partial_{\nu}w^{i,k}\phi_i-w^{i,k}\partial_{\nu}\phi_i)+\sum_i \int_{\tilde \Omega_k}(\Delta \phi^i+e^{U_i}\phi_i)w_i^k. $$
The second term is $0$ according to the equation for $\phi^i$. We claim that the first term is $O(\epsilon_k^{3})$ if $\phi_i$ is of the following form:
\begin{equation}\label{phif2}
\phi^i=\frac 1{r^2}(d_i\cos (2\theta)+q_i\sin(2\theta))+O(1/r^3).
\end{equation}
Indeed, using (\ref{phif2}) and $\nabla w_i^k=O(\epsilon_k^2)$ on $\partial \tilde \Omega_k$, we see that
$$\int_{\partial \tilde \Omega_k}\partial_{\nu}w^{i,k}\phi_i=O(\epsilon_k^3). $$
On the other hand, by $|w_i^k(y)|\le C\epsilon_k$ on $\partial \tilde \Omega_k$ we have
$$\int_{\partial \tilde \Omega_k}w^{i,k}\partial_{\nu}\phi_i=O(\epsilon_k^{3}). $$
Therefore the left hand side is $O(\epsilon_k^3)$. Now we estimate terms on the right hand side.
Using the equation for $\phi^i$ and (\ref{phif2}) we have
 \begin{align}\label{phif2a}
&\int_{\tilde \Omega_k}r^2e^{U^k_i}\phi_i\cos (2\theta)
=-\int_{\tilde \Omega_k}r^2\cos(2\theta)\Delta \phi^i \nonumber\\
=&-\int_{\partial \tilde \Omega_k}(\partial_{\nu}\phi^i-\frac{2\phi^i}r)r^2\cos(2\theta)dS
=-\int_0^{2\pi}(\partial_r\phi^i-\frac{2\phi^i}r)r^3\cos(2\theta)d\theta \nonumber\\
=&4\pi d_i+O(\epsilon_k).
\end{align}
Similarly
\begin{equation}\label{phif2b}
\int_{\tilde \Omega_k}r^2e^{U^k_i}\phi_i\sin (2\theta)=4\pi q_i+O(\epsilon_k).
\end{equation}

Using the vanishing rate of $\nabla h^k_i(0)$ and (\ref{phif2}) we have
\begin{equation}\label{phif2c}
\epsilon_k \int_{\tilde \Omega_k}e^{U^k_i}\phi_i\partial_1h^k_i(0)y_1 dy =O(\epsilon_k^3),
\end{equation}
and
\begin{equation}\label{phif2d}
\epsilon_k \int_{\tilde \Omega_k}e^{U^k_i}\phi_i\partial_2h^k_i(0)y_2 dy =O(\epsilon_k^3).
\end{equation}
Finally by (\ref{phif2}) we have
\begin{align}\label{phif2e}
&\int_{\tilde \Omega_k}r^2e^{U^k_i}\phi_idy=\int_{\tilde \Omega_k}(-\Delta \phi^i)r^2 \\
=&\int_{\partial \tilde \Omega_k}(-\partial_{\nu}\phi^ir^2+2\phi^i r)dS-4\int_{\tilde \Omega_k}\phi^i \nonumber \\
=&-4\int_{\mathbb R^2}\phi^i +O(\epsilon_k) \nonumber
\end{align}
where in the last line above we used $\int_{\tilde \Omega_k}\phi^i=\int_{\mathbb R^2}\phi^i+O(\epsilon_k)$ because for $r$ large
(\ref{phif2}) holds, the integration of the leading term is $0$ and the integration of the error term gives $O(\epsilon_k)$.

Combing (\ref{phif2}),(\ref{phif2a}),(\ref{phif2b}),(\ref{phif2c}),(\ref{phif2d}) and (\ref{phif2e}) we have,
\begin{equation}\label{imeq}
\sum_{i=1}^n\bigg (\Delta h^k_i(0)(\int_{\mathbb R^2}\phi^i)-d_i\pi(\partial_{11}h^k_i(0)-\partial_{22}h^k_i(0))-2q_i\pi \partial_{12}h^k_i(0)\bigg )=O(\epsilon_k)
\end{equation}

Now we defined $2n-2$ families of solutions to the linearized Toda system: For $m=1,...,n-1$,
$$\Phi_m=\left(\begin{array}{c}\phi^1_m\\
\vdots \\
\phi^n_m
\end{array}
\right )=\frac{-\partial \bf U^k}{\partial \alpha^k_{m+1,2}}=
\left(\begin{array}{c}
\frac{-\partial U^{1,k}}{\partial \alpha^k_{m+1,2}}\\
\vdots \\
\frac{-\partial U^{n,k}}{\partial \alpha^k_{m+1,2}}
\end{array}
\right )=\left(\begin{array}{c}O(r^{-3})\\
\vdots\\
\frac{m(m+1)}{r^2}\cos(2\theta)+O(r^{-3})\\
-\frac{m(m+1)}{r^2}\cos(2\theta)+O(r^{-3})\\
\vdots \\
O(r^{-3})
\end{array}
\right )
$$
and for $m=n,...,2n-2$
$$\Phi_m=\frac{-\partial \bf U^k}{\partial \beta^k_{m+1,2}}=\left(\begin{array}{c}O(r^{-3})\\
\vdots\\
\frac{m(m+1)}{r^2}\sin(2\theta)+O(r^{-3})\\
-\frac{m(m+1)}{r^2}\sin(2\theta)+O(r^{-3})\\
\vdots \\
O(r^{-3})
\end{array}
\right )
$$
where major entries appear in row $m$ and row $m+1$ in both vectors. Here we also recall that $c^k_{n+2-m,n-m}=\alpha_{m,2}^k+\sqrt{-1}\beta_{m,2}^k$,
$c_{n+1-m}^k=\alpha_{m+1,2}^l+\sqrt{-1}\beta_{m+1,2}^2$.

Using (\ref{ualpha}) and (\ref{ubeta}) in (\ref{imeq}) we see that (\ref{2van}) holds. Theorem \ref{localthm} is proved. $\Box$

\medskip

As an immediate consequence of Theorem \ref{localthm} we have

\begin{cor}\label{cor-2} With the same assumptions of Theorem \ref{localthm}, we have
\begin{equation}\label{log-term}
\int_{B_{\sigma}}h_ie^{u_i^k}dx=4\pi i(n+1-i)+\pi \Delta(\log h_i)(0)e^{a_{i,k}}\epsilon_k^2\log \frac{1}{\epsilon_k}+O(\epsilon_k^2)
\end{equation}
where $a_{i,k}$ is the leading coefficient of $U_{i,k}$:
$$U_{i,k}(x)=-4\log |x|+a_{i,k}+O(r^{-1}), \quad \mbox{ for }\quad |x|>1. $$
\end{cor}

\noindent{\bf Proof of Corollary \ref{cor-2}:}
Recall that $h_i^k(x)=\frac{h_i(x)}{h_i(0)}e^{\phi_i^k(x)}$.
\begin{align*}
&\int_{B_{\sigma}}h_ie^{u_i^k}dx=\int_{B_1}h_i^ke^{u_i^k+\log h_i(0)-\phi_i^k}dx\\
=&\int_{B(0,\epsilon_k^{-1}\sigma)}h_i^k(\epsilon_k y)e^{\tilde v_i^k(y)}dy\\
=&\int_{B(0,\epsilon_k^{-1})}\bigg (1+\epsilon_k \nabla h_i^k(0)y \\
&+\epsilon_k^2\big (\partial_{11}h_i^k(0)(y_1^2-\frac{r^2}2)
+\partial_{22}h_i^k(0)(y_1^2-\frac{r^2}2)+2\partial_{12}h_i^k(0)y_1y_2\big ) \bigg )e^{\tilde v_i^k(y)}\\
&+\frac 12\Delta h_i^k(0)r^2\epsilon_k^2e^{\tilde v_i^k(y)}dy+O(\epsilon_k^2).
\end{align*}
By Theorem \ref{localthm} we have
$\nabla h_i^k(0)=O(\epsilon_k)$,
$\Delta h_i^k(0)=\Delta (\log h_i)(0)$, and
$$e^{\tilde v_i^k(y)}=e^{U_i^k(y)}(1+O(\epsilon_k^2)(1+|y|))\quad \mbox{ for }|y|>1 .$$
Then elementary computation leads to  (\ref{log-term}). Corollary \ref{cor-2} is established. $\Box$

\section{Proof of the main theorems}

Without loss of generality we assume
\begin{equation}\label{nliz}
\int_M h_ie^{u_i^k}dV_g=1,\quad i=1,...,n.
\end{equation}
because otherwise we just use $\tilde u_i^k=u_i^k-\log \int_M h_ie^{u_i^k}dV_g$. For simplicity we just assume that $u_i^k$ satisfies (\ref{nliz}).

Recall that $p_1,..,p_L$ are distinct blowup points and $p_t^k$ satisfies
$$\max_iu_i^k(p_t^k)=\max_i \max_{B(p_t^k,\delta)}u_i^k$$
 where $B(p_t^k,\delta)$ ($t=1,..,L$) are mutually disjoint balls.
 Around each blowup point, say, $p_1^k$, we use the local coordinate, then $ds^2=e^{\psi(y_{p_1^k})}(dy_1^2+dy_2^2)$ where
$$\nabla \psi(0)=0,\quad \psi(0)=0,\quad \Delta \psi=-2K e^{\phi} $$
where $K$ is the Gauss curvature.

In local coordinates around $p_1^k$, we write the system (\ref{main-uk}) as
\begin{equation}\label{uiklc}
\Delta u_i^k+\sum_{j=1}^n a_{ij}\rho_j^ke^{\psi}(h_je^{u_j^k}-1)=0,\quad B_{\delta}.
\end{equation}
where (\ref{nliz}) is used.

Let $f_i$ satisfy (for simplicity we omit $k$ in this notation)
\begin{equation}\label{fik}
\left\{\begin{array}{ll}
\Delta f_i-\sum_j a_{ij}\rho_j^ke^{\psi}=0, \quad B_{\delta}, \\
f_i=0,\quad \mbox{ on }\quad \partial B_{\delta}.
\end{array}
\right.
\end{equation}
Let $\tilde u_i^k=u_i^k-f_i$, then we have
\begin{equation}\label{tuik}
\Delta \tilde u_i^k+\sum_{j=1}^n a_{ij}\tilde h_j^k e^{\tilde u_j^k}=0,\quad B_{\delta}
\end{equation}
where
\begin{equation}\label{th}
\tilde h_i^k=\rho_i^k e^{\psi}h_ie^{f_i}.
\end{equation}
The boundary oscillation of $\tilde u_i^k$ is finite, because by the Green's representation of $u_i^k$ we have:
\begin{equation}\label{gr1}
u_i^k(x)=\bar u_i^k+\int_M G(x,\eta)\sum_{j=1}^n a_{ij}\rho_j^k h_j e^{u_j^k}dV_g.
\end{equation}
Since $u_i^k$ is bounded above away from bubbling areas:
$$u_i^k(x)\le C \quad \forall x\in M\setminus \cup_{t=1}^LB(p_t^k,\delta), $$
it is easy to see that the oscillation of $u_i^k$ on $M\setminus \cup_{t=1}^LB(p_t^k,\delta)$ is finite.

Let $\tilde M_k=\max_{x\in B(p_t,\delta)}\max_i u_i^k$ and $\tilde \epsilon_k=e^{-\frac 12 \tilde M_k}$,
$$\tilde v_i^k(y)=\tilde u_i^k(\tilde \epsilon_k y+p_1^k)+2\log \tilde \epsilon_k. $$
By Theorem  \ref{localthm} there exists a sequence of $\tilde U_i^k$ such that
$$\tilde U_i^k(y)=-4\log (1+|y|)+O(1),\quad y\in \mathbb R^2 $$
and
$$|\tilde v_i^k(y)-\tilde U_i^k(y)|\le C(\delta)\tilde \epsilon_k(1+|y|),\quad |y|\le \tilde \epsilon_k^{-1}\delta. $$
Thus
$$u_i^k(x)=-2\log \tilde \epsilon_k +O(1), \mbox{ for } \quad x\in B(p_1^k,\delta)\setminus B(p_1^k, \delta/2). $$
Applying this argument to each bubbling area we see that
$$\max_i \max_{B(p_t,\delta)}u_i^k=\max_i \max_{B(p_l^k,\delta)}u_i^k+O(1),\quad \forall t\neq l. $$
Thus we shall just use
$$M_k=\max_i \max_{x\in M}u_i^k, \quad \mbox{ and } \quad \epsilon_k=e^{-\frac 12 M_k} $$
from now on. We also have derived from the argument above that
$$\bar u_i^k=-M_k+O(1)$$
and
\begin{equation}\label{uaway}
u_i^k(x)=-M_k+O(1),\quad \mbox{ for } \, i=1,...,n, \quad x\in M\setminus \cup_{t=1}^L B(p_t^k,\delta).
\end{equation}

Next we evaluate $u_i^k$ on $\partial B(p_1^k,\delta)$, using (\ref{uaway}) and (\ref{gr1}) we have
\begin{align}
u_i^k(x)&=\bar u_i^k+\int_M G(x,\eta)(\sum_j a_{ij}\rho_j^kh_je^{u_j^k})dV_g \nonumber\\
=&\bar u_i^k+\sum_{t=1}^L \int_{B(p_t,\delta)}G(x,\eta)(\sum_j a_{ij}\rho_j^k h_je^{u_j^k})dV_g+O(\epsilon_k^2). \label{gr2}
\end{align}
From (\ref{th}) we see that
\begin{equation}\label{thik}
dV_g=e^{\psi}dx, \quad \rho_i^kh_ie^{\psi}e^{u_i^k}=\tilde h_i^ke^{\tilde u_i^k}.
\end{equation}
Thus by Theorem \ref{localthm}
$$\int_{B_{\delta}}\sum_j a_{ij}\tilde h_j^ke^{\tilde u_j^k}=\int_{B(0,\epsilon_k^{-1}\delta)}\sum_j a_{ij}\tilde h_j^k(\epsilon_k y)
e^{\tilde v_j^k(y)}dy=8\pi +O(\epsilon_k). $$
Moreover by similar computation we have
$$\int_{B(p_1^k,\delta)}(G(x,p_1^k)-G(x,\eta))\sum_j a_{ij} \rho_j^k h_j e^{u_j^k}dV_g=O(\epsilon_k). $$
Thus
\begin{align}
&u_i^k(x)=\bar u_i^k+8\pi G(x,p_1^k)+8\pi \sum_{t=2}^LG(x,p_t^k)+O(\epsilon_k)\nonumber\\
&=\bar u_i^k-4\log |x-p_1^k|+8\pi \gamma(x,p_1^k)+8\pi \sum_{t=2}^LG(x,p_t^k)+O(\epsilon_k) \label{gr3}
\end{align}
for $x$ close to $p_1^k$.

From (\ref{gr3}) we can determine which harmonic function $\phi_i^k$ that makes $u_i^k-\phi_i^k$ a constant on $\partial B(p_1^k,\delta)$. Before we determine $\phi_i^k$ we first observe that
\begin{equation}\label{rho-crude}
\rho_i^k=4\pi i(n+1-i)L+O(\epsilon_k).
\end{equation}
Indeed,
\begin{align}
&\rho_i^k=\int_M \rho_i^k h_i e^{u_i^k} dV_g \nonumber\\
=&\sum_{t=1}^L \int_{B(p_t^k,\delta)}\rho_i^k h_i e^{u_i^k} dV_g+O(\epsilon_k^2)\nonumber \\
=&\sum_{t=1}^L \int_{B(p_t^k,\delta)} \tilde h_i^k e^{\tilde u_i^k}dx+O(\epsilon_k^2). \label{tem-1}
\end{align}
By Theorem \ref{localthm}
$$\int_{B(p_t^k,\delta)}\sum_j a_{ij} \tilde h_j^k e^{\tilde u_j^k}dx=8\pi +O(\epsilon_k), $$
which implies
\begin{equation}\label{tem-2}
\int_{B(p_t^k,\delta)} \tilde h_i^k e^{\tilde u_i^k}dx=4\pi i(n+1-i) +O(\epsilon_k),
\end{equation}
thus (\ref{rho-crude}) is verified by (\ref{tem-1}) and (\ref{tem-2}).

Let $\phi_i^k$ be a sequence of harmonic functions that annihilates the oscillation of $u_i^k$ on $\partial B_{\delta}$:
$$\left\{\begin{array}{ll}
-\Delta \phi_i^k=0,\quad \mbox{ in } \quad B_{\delta}, \\
\phi_i^k(x)=u_i^k(x)-\frac{1}{2\pi \delta}\int_{\partial B_{\delta}}u_i^kdS, \quad \mbox{ for } x\in \partial B_{\delta}.
\end{array}
\right.
$$
Then we claim that
\begin{align}\label{phi-near-1}
\phi_i^k(x)=&8\pi(\gamma(x,p_1^k)+\sum_{t=2}^LG(x,p_t^k))-f_i\\
&-8\pi(\gamma(p_1^k,p_1^k)+\sum_{t=2}^L G(p_1^k,p_t^k))+f_i(0)+O(\epsilon_k).\nonumber
\end{align}

Indeed, from (\ref{greenf}) we see that
$$-\Delta G(x,p_1^k)=\delta_{p_1^k}-e^{\psi}, $$
comparing with (\ref{fik}) we have
$$\Delta \bigg (8\pi\big (\gamma(x,p_1^k)+\sum_{t=2}^LG(x,p_t^k)\big )-\frac{\sum_ja_{ij}\rho_j^k}{8\pi L}f_i\bigg )=0. $$
In (\ref{gr3}) we see that $-4\log |x-p_1^k|$ has no oscillation on $\partial B(p_1^k,\delta)$, neither has $f_i$. Thus (\ref{phi-near-1}) is verified.
Correspondingly we also have
\begin{equation}\label{near-1d}
\nabla \phi_i^k(0)=8\pi \bigg (\nabla_1 \gamma(p_1^k,p_1^k)+\sum_{t=2}^L\nabla_1G(p_1^k,p_t^k)\bigg )-\nabla f_i^k(p_k)+O(\epsilon_k).
\end{equation}
where $\nabla_1$ means differentiation with respect to the first component.
Going back to (\ref{tuik}) and applying Theorem \ref{localthm}, we have
$$\frac{\nabla h_i(p_1^k)}{h_i(p_1^k)}+8\pi\bigg (\nabla_1 \gamma(p_1^k,p_1^k)+\sum_{t=2}^L\nabla_1 G(p_1^k,p_t^k)\bigg )=O(\epsilon_k).
$$
From (\ref{thik}) we have
$$\log \tilde h_i^k(x)=\log \rho_i^k+\psi+\log h_i+f_i, $$
and
$$\Delta (\log \tilde h_i^k(0))=-2K(p_1^k)+\Delta (\log h_i)(p_1^k)+8\pi L+O(\epsilon_k). $$
Let
\begin{align*}
H_i^k&=\log \tilde h_i^k+\phi_i^k\\
&=\log \rho_i^k+\psi+\log h_i+8\pi(\gamma(x,p_1^k)+\sum_{t=2}^LG(x,p_t^k))\\
&-8\pi\bigg(\gamma(p_1^k,p_1^k)+\sum_{t=2}^LG(p_1^k,p_t^k)\bigg )+f_i(0)+O(\epsilon_k).
\end{align*}
$H_i^k$ corresponds to $h_i^k$ in the statement of Theorem \ref{localthm}. Clearly
\begin{align*}
\Delta H_i^k(0)&=-2K(p_1^k)+\Delta (\log h_i)(p_1^k)+8\pi L+O(\epsilon_k)\\
(\partial_{11}-\partial_{12})(H_l^k-H_{l-1}^k)(0)&=(\partial_{11}-\partial_{12})(\log h_{l-1}-\log h_l)(p_1^k)\\
\partial_{12}(H_l^k-H_{l-1}^k)(0)&=\partial_{12}(\log h_l-\log h_{l-1})(p_1^k).
\end{align*}
Using the expressions above we obtain (\ref{2ndf1}) and (\ref{2ndf2})
from the vanishing estimate for second frequencies in Theorem \ref{localthm} around $p_1^k$. The corresponding $T_l^k$ and $\tilde T_l^k$ are defined by the approximating function $U_{i,1}^k$, which satisfies
$$|u_i^k(\epsilon_ky+p_1^k)+\log h_i(p_1^k)+\phi_i^k(\epsilon_k y)-U_{i,1}^k(y)|\le C\epsilon_k(1+|y|),\quad |y|\le \epsilon_k^{-1}\delta. $$
For $t=1,...,L$, we have $U_{i,t}^k$ as the approximating global functions in $B(p_t^k,\delta)$ with parameters $\lambda_{i,t}^k$ for $i=0,..,n$, $c_{ij,t}^k$.  For $l=2,...,n$
\begin{align}\label{Tlk}
&T_{l-1,k}^t=\int_{\mathbb R^2}(-\frac{\partial U_t^{l-1,k}}{\partial \alpha_{l,2,t}^k}), \qquad
\tilde T_{l-1,k}^t=\int_{\mathbb R^2}(-\frac{\partial U_t^{l-1,k}}{\partial \beta_{l,2,t}^k}), \\
&T_{l,k}^t=\int_{\mathbb R^2}(-\frac{\partial U_t^{l,k}}{\partial \alpha_{l,2,t}^k}),\qquad
\tilde T_{l,k}^t=\int_{\mathbb R^2}(-\frac{\partial U_t^{l,k}}{\partial \beta_{l,2,t}^k}),\quad l=2,...,n,\nonumber \\
&U^{i,k}_t=\sum_{j=1}^n a^{ij}U_{j,t}^k, \quad i=1,...,n \nonumber \\
&c_{n+2-i,n-i,t}^k=\alpha_{i,2,t}^k+\sqrt{-1}\beta_{i,2,t}^k,\quad i=2,...n. \nonumber
\end{align}
From the expansion of $U_t^{i,k}$ in (\ref{euk2}) it is easy to see that all the integrals are finite and are $O(\epsilon_k)$ different from the integration on $B(0,\delta \epsilon_k^{-1})$.

Finally Theorem \ref{rho-1} is direct consequence of Corollary \ref{cor-2}.
Main theorems (Theorem \ref{rho-1}, Theorem \ref{loc-bu-1} and Theorem \ref{loc-bu-2}) are established. $\Box$

\begin{rem} It can be verified that when $n=2$ and $l=2$, the $\partial_z^2$ condition (\ref{2ndf1})\&(\ref{2ndf2}) is the same as the (14)\&(15) in \cite{lin-wei-zhao} if the error in (14)\&(15) of \cite{lin-wei-zhao} is replaced by $O(\epsilon_k)$.
\end{rem}

\end{document}